\input amstex

\loadeufm
\loadmsbm
\loadeufm

\documentstyle{amsppt}
\input amstex
\catcode `\@=11
\def\logo@{}
\catcode `\@=11
\magnification \magstep1
\NoRunningHeads
\NoBlackBoxes
\TagsOnLeft

\def \={\ = \ }
\def \+{\ +\ }
\def \-{\ - \ }

\def \b|{\big |}

\def \g1{\Gamma_1}

\def \nfp{\demo\nofrills{Proof:\usualspace\usualspace }}

\def\rarr#1#2{\smash{\mathop{\hbox to .5in{\rightarrowfill}}
 	 \limits^{\scriptstyle#1}_{\scriptstyle#2}}}

\def\larr#1#2{\smash{\mathop{\hbox to .5in{\leftarrowfill}}
	  \limits^{\scriptstyle#1}_{\scriptstyle#2}}}

\def\swarr#1#2 {\llap{$\scriptstyle #1$}  \swarrow
  	\vcenter to .5in{}\rlap{$\scriptstyle #2$}}

\topmatter
\title Automorphisms of groups and a higher rank JSJ decomposition I: RAAGs and a higher rank Makanin-Razborov diagram 
\endtitle
\author
\centerline{ 
Z. Sela${}^{1,2}$}
\endauthor

\dedicatory{Dedicated to Thomas Delzant on his 60th birthday}
\enddedicatory

\footnote""{${}^1$Hebrew University, Jerusalem 91904, Israel.}
\footnote""{${}^2$Partially supported by an Israel academy of sciences fellowship.} 
\abstract\nofrills{}
The JSJ decomposition encodes the automorphisms and the virtually cyclic splittings
of a hyperbolic group. For general finitely presented groups, the JSJ decomposition
encodes only their splittings. 

In this sequence of papers we study the automorphisms 
of a hierarchically hyperbolic group that satisfies some weak acylindricity conditions.
To study these automorphisms we construct an object that can be viewed as a higher rank JSJ
decomposition.

In the first paper we demonstrate our construction in the case of a right angled Artin group. 
For studying automorphisms of a general HHG we construct what we view as a higher rank Makanin-Razborov diagram, which 
is the first step in the construction of the higher rank JSJ.
\endabstract
\endtopmatter

\document

\baselineskip 12pt

The (canonical) JSJ decomposition of a  torsion-free hyperbolic group was originally constructed to generalize the solution
of the isomorphism problem from rigid torsion-free hyperbolic groups to all torsion-free hyperbolic groups. Such a generalization
required an understanding and techniques to handle both automorphisms and splittings of such groups, and both are encoded by the
JSJ decomposition ([Se3],[Le]). 

The construction of the JSJ was later generalized to general finitely presented groups (see [Gu-Le]). In this general setting, the JSJ encodes
all the splittings of a f.p.\ group over a given family of subgroups
(in a rather subtle way), 
but it is far from encoding the automorphism group nor
the dynamics of individual automorphisms.

In this sequence of papers we use some of the JSJ concepts to study automorphisms of hierarchically hyperbolic groups. Hierarchically hyperbolic groups and spaces
were defined by Behrstock, Hagen and Sisto [BHS1]. The definition axiomatizes the hierarchical structure of the mapping class groups, that was defined
and studied in the work of Masur and Minsky [Ma-Mi]. Automorphisms of families of HHG were studied earlier by  
Fioravanti [Fi], and by Casals-Ruiz, Hagen and Kazachkov ([Ca-Ka],[CHK]). 

To study the automorphism group and the dynamics of individual automorphisms of an HHG, we look at the action of some characteristic
finite index subgroup of the HHG
on the domains that are part of the HHG structure. We further associate a virtually abelian decomposition of some quotient of the 
finite index characteristic subgroup of the HHG with each orbit
of  domains under the action of the finite index subgroup. The (finite) collection of virtually abelian decompositions that we construct can
be viewed as a higher rank JSJ decomposition of the HHG. 

The finite collection of decompositions encode the dynamics of individual automorphisms, and can
be used to study the algebraic structure of the (outer) automorphism group of the HHG. To construct the virtually abelian decompositions, we borrow techniques that
were previously used to study sets of solutions to system of equations (varieties) over certain families of groups (e.g.\ [Se1], [Ja-Se], [Re-We]  and [Gr-Hu]),
together with basic tools that were used in the  study of the first order theory of the free group (that appear in [Se2] and [Se4]).

To be able to apply these techniques and constructions we require that the action of the set stabilizer of each domain in the HHG  on the domain that it stabilizes
is $weakly$ $acylindrical$ (definition 2.1). This is a coarse form of an acylindrical action of a set stabilizer modulo the pointwise stabilizer of the space.
We further require that the HHG is colorable, i.e., that it has a finite index subgroup, 
for which the domains in each orbit of the finite index subgroup are 
pairwise transverse. This is known to be true for the mapping class groups [BBF], and in many other cases (see [Ha-Pe] and [DMS]). 

In the first section of the paper we motivate and demonstrate our approach by examining automorphisms of right angled Artin groups, based on works
of Charney-Crisp-Vogtmann ([CCV] and [CV]), and Duncan-Kazachkov-Remeslennikov [DKR]. With each RAAG we associate a finite collection of graphs of groups
and simplicial trees that are associated with them, on which
the RAAG acts weakly acylindrically. Furthermore, the automorphism group of the RAAG preserves free products  that are associated with the graphs of groups
(proposition 1.5),
so it is not difficult to read from these graphs of groups and their  Bass-Serre trees the structure of the higher rank JSJ decomposition 
(we do not present these last constructions in this paper).

In sections 2 and 3 we construct what we view as a higher rank Makanin-Razborov diagram that is associated with a general (colorable) 
HHG that satisfies our weak acylindricity
assumption (Theorems 3.3 and  3.4). The construction of the higher rank diagram is not canonical. It is based on a compactness argument,
 and follows the steps of the construction of such diagrams in [Ja-Se], [Re-We] and 
[Gr-Hu]. 

The diagram that we construct is not canonical but it is universal. Every automorphism of the HHG factors through at least one of its collections of (cover)
resolutions (see definition 2.6 for the precise definition of an automorphism that factors through a resolution). 
The existence of the higher rank diagram is the basis for the construction of the higher rank JSJ in the next papers. In section 2 we construct
a higher rank diagram in case the HHG is a product of hyperbolic spaces (Theorem 2.7), and in section 3 we generalize the construction to HHG that satisfy
our assumptions (theorems 3.3 and  3.4).

The whole approach that we adapt in this sequence of papers aims to demonstrate that techniques that were originally
developed to analyze homomorphisms and varieties over groups,
and more generally to study first order formulas, can be used to study automorphisms, that are in general transcadental and are certainly not (first order) definable.
Still the arguments that appear in this (first) paper, can analyze homomorphisms into HHGs as well. 

In the last section we state an analogue of theorem 3.4
for the structure of homomorphisms of a f.g.\ group into a (colorable)  HHG that satisfies our weak acylindricity assumption. Finally, we apply our theorem
to get another generalization of theorem 3.4, that describes all the homomorphisms of a given f.g.\ group into all the (colorable) HHGs that satisfy 
some uniform HHS and weak acylindricity assumptions. In
some sense we view that generalization
as an analogue of Thurston's (bounded image) 
theorem on the geometric structure of discrete faithful representations of a f.g.\ group into a rank 1 Lie group, and its connection to the JSJ decomposition
(see [Mo], [Se3] and [Ka]),
to representations into a uniform
family of (colorable) HHGs.   

\smallskip
This whole project started as an attempt to answer a question of Eliyahu Rips on the structure of the automorphism groups of cubulated groups, and
their connection to automorphisms of low dimensional manifolds, along the line of the JSJ decomposition of a hyperbolic group. It is also a late
answer to a question of Ruth Charney on the possibility to encode the automorphisms of a RAAG by a JSJ decomposition. I am indebted to both of them.
Eran Segalis  found a mistake in an earlier formulation of proposition 1.5, and the referees reports were full of insight and excellent comments and advices.
I am grateful to them and to the editor.

\vglue 1.5pc
\centerline{\bf{\S1. The flags hypergraph and automorphisms of a RAAG}}
\medskip
Let $\Gamma$ be a finite simple graph, and let $A_{\Gamma}$ be the right-angled Artin group that
is associated with $\Gamma$. 
$A_{\Gamma}$ is freely indecomposable if and only if $\Gamma$ is connected. The decomposition that we introduce 
does not give any new insight in studying freely decomposable groups. Hence, for the rest of this section we will
assume that the graph $\Gamma$ is connected, and to avoid trivialities we also assume that $\Gamma$ contains at least 3 vertices.

In [CCV] and [CV] the outer automorphism group of a  RAAG $A_{\Gamma}$ is studied by mapping a finite index subgroup of it
into the direct sum of outer automorphisms of maximal $joins$. Joins $U*V$ are complete bi-partite graphs on the set
of vertices $U$ and $V$.  This map gives a bound on the virtual cohomological dimension of $Out(A_{\Gamma})$, and
enabled the authors to construct an outer space on which a finite index subgroup of $Out(A_{\Gamma})$ acts.

We note that our results for RAAGs are based on the works of Charney-Vogtmann [CV] and Duncan-Kazachkov-Remeslennikov [DKR].
We mainly present their results in a different way, that demonstrates our approach to the construction of the higher rank JSJ decomposition for
more general HHGs in the sequel.

We are interested in the specific description of the group $Out(A_{\Gamma})$. We will associate a (canonical)
$flags$ $hypergraph$  with $\Gamma$, that is related to the set of maximal joins that is
studied in [CCV] and [CV], but our point of view is somewhat different. We look at RAAGs as an example, or a motivation, for the general
definition of a higher rank JSJ decomposition. And we view the flags hypergraph as encoding the higher
rank JSJ decomposition of a  RAAG, that will give us the description of $Out(A_{\Gamma})$. 

Following [CCV] we let $v \in \Gamma$ be a vertex, and denote $lk(v)$ the $link$ of $v$, and  $st(v)$ the $star$ of $v$.
Following [CV], on the set of vertices of $\Gamma$ we define a partial order.
Given $u,v \in \Gamma$, we say that $u \leq v$ if $lk(u) \subset st(v)$. 
In ([CV], lemma 2.2) it is proved that this
is a partial order. With the partial order one can naturally associate an equivalence classes of vertices, that we denote $[u]$.

Note that if $[u]=[v]$, and there is no edge between $u$ and $v$ in $\Gamma$, then $lk(u)=lk(v)$. Hence, there are no edges between the vertices
in $[u]$, the vertices in $[u]$ form an anti-clique, and the generators in $[u]$ generate a free subgroup in $A_{\Gamma}$. 
If $[u]=[v]$ and there is an edge between $u$ and $v$ in $\Gamma$, then $st(u)=st(v)$, and there exists an edge between any two vertices in $[u]$. In
this case the vertices in $[u]$ form a clique, and $A_{[u]}$ is free abelian.

In the sequel we will use the notation,  $lk([u])=lk(u) \setminus \, [u]$ and $st([u])=lk(u) \cup [u]$, 
for any of the vertices $u \in [u]$.
The results of [CV] are obtained by studying the maximal equivalence classes w.r.t.\ the above partial order. To
get the hypergraph we are looking for, and obtain the structure of $Out(A_{\Gamma})$, we need to look at all 
the equivalence classes and not only the maximal ones.

From the partially ordered set of equivalence classes of vertices in $\Gamma$, we can construct the flags hypergraph,
$\Delta_{\Gamma}$,
that is associated with a  RAAG $A_{\Gamma}$.

\vglue 1.5pc
\proclaim{Definition 1.1} 
The vertices in the hypergraph $\Delta_{\Gamma}$ are  the vertices of the graph $\Gamma$.
We start the 
construction of $\Delta_{\Gamma}$ with the collection of the maximal equivalence classes with respect to the partial order
that was defined on the vertices of $\Gamma$. We view each maximal equivalence class as a hyperedge of level 1, and each vertex
in a maximal equivalence class as a vertex of level 1.

At step $k$, we look at all the classes $[v]$ that are maximal after we took out all the classes
of level up to $k-1$. With each such class $[v]$ we associate a hyperedge of level $k$. This hyperedge contains all the hyperedges
of level less than $k$, such that all the classes in these hyperedges are bigger than $[v]$, and all the vertices in $[v]$. The vertices in
$[v]$ are defined to be of level $k$ as well. Since $\Gamma$ is finite, the construction of the hypergraph $\Delta_{\Gamma}$ terminates
after finitely many steps.

We say that a hyperedge $E_k$ of level $k \geq 2$ is $centerless$ if its highest level vertices centralize vertices in $E_k$ only if they are
 in their class. In that
case: $A_{E_k}=B*A_{[v]}$, where $[v]$ is the equivalence class of the highest level vertices in $E_k$, and $B$ is the subgroup generated by the complement
of $[v]$ in $E_k$.

A hyperedge $E_k$ of level $k \geq 2$
that is not centerless must have a non-trivial center. If $[v]$ is the class of the highest level vertices in $E_k$, then the center of $A_{E_k}$ 
is a free abelian group that is generated by all the vertices in the complement of $[v]$ in $E_k$ that commute with the vertices in $[v]$.     
In that case $A_{E_k}=Ab \oplus(B * A_{[v]})$, where $Ab$ is in the center of $A_{E_k}$, and $B$ is the (possibly trivial) subgroup that is generated by all the generators
in $E_k$ that are not in $[v]$ and not in  $Ab$. Note that in case $B$ is trivial and $A_{[v]}$ is free abelian, $A_{E_k}$ is abelian. Otherwise,
$Ab$ is the center of $A_{E_k}$. We say that such a hyperedge is a hyperedge $with$ $center$. 
\endproclaim

By work of M. Laurence [La], building on work of Servatius [Ser], the automorphism group of a RAAG, $Aut(A_{\Gamma})$, is generated by:
\roster
\item"{(1)}" Inner automorphisms.

\item"{(2)}" Symmetries induced by symmetries of the graph $\Gamma$, that permute the standard generators.

\item"{(3)}" Inversions that send a standard generator to its inverse.

\item"{(4)}" Transvections. Whenever $v \leq w$, a transvection sends $v$ to $vw$.

\item"{(5)}" Partial conjugations. If $st(v)$ separates $\Gamma$, a partial conjugation conjugates all the generators
in a connected component of $\Gamma \setminus \, st(v)$ by $v$.
\endroster

With the flags hypergraph $\Delta_{\Gamma}$ we can associate finitely many  groups of  automorphisms of $A_{\Gamma}$. 

\vglue 1.5pc
\proclaim{Definition 1.2} Let $\Gamma$ be a finite graph. If $\Gamma$ is disconnected, $A_{\Gamma}$ is freely decomposable. In that case,
 we associate a group of automorphisms with each connected component, and then add the automorphisms of the corresponding free product. Hence,
we assume that $\Gamma$ is connected.

Let $\Delta_{\Gamma}$ be its flags 
hypergraph. We add a group of automorphisms for each hyperedge in $\Delta_{\Gamma}$, according to their natural grading - their levels.
On the vertices of $\Delta_{\Gamma}$ there is a natural grading - their levels. 
We start with hyperedges of level $1$ in $\Delta_{\Gamma}$.

Let $[u]$ be  an equivalence class of vertices of level $1$ in $\Gamma$ (note that $[u]$ is a hyperedge of level 1 in $\Delta_{\Gamma}$). By [CV] $A_{[u]}$
is either free or free abelian with a free (abelian) basis, the generators that are associated with the vertices in [u]. With $[u]$
we associate a group of automorphisms of $A_{\Gamma}$ that depends if [u] is free or free abelian. 
The group of automorphisms that we add for hyperedges of level 1 are generated by:
\roster
\item"{(1)}" In case $A_{[u]}$ is free, the automorphism group of the free group that is generated by the vertices in $[u]$, $Aut(A_{[u]})$. In case $A_{[u]}$ is
free abelian, its automorphism group (isomorphic to $GL(n,Z)$). Each such automorphism extends to an
automorphism of $A_{\Gamma}$ by defining it to be the identity on all the generators that are not in $[u]$.

\item"{(2)}" We look at the complement of the star of $[u]$ in $\Gamma$: $\Gamma \setminus st([u])$. Let $\Gamma^{[u]}_1,\ldots,\Gamma^{[u]}_m$ be 
the connected components of that complement that are not single vertices. We join to automorphisms of type (1) the automorphisms of $A_{\Gamma}$ that are obtained by
conjugating each of the subgroups, $A_{\Gamma^{[u]}_i}$, $i=1,\ldots,m$, by elements from $A_{[u]}$. Note that these are particular automorphisms of the free
product:
$$A_{[u]}*A_{\Gamma^{[u]}_1}*\ldots*A_{\Gamma^{[u]}_m}$$
that extend naturally to automorphisms of $A_{\Gamma}$ by defining them to be the identity on the generators that are connected to $[u]$, and on generators that 
are roots that are connected to vertices in $lk([u])$.
\endroster

We continue to the higher level hyperedges. 
Let $E_k$ be a centerless hyperedge of level $k$ in $\Delta_{\Gamma}$. 
The hyperedge $E_k$ contains vertices of lower level and an equivalence class $[v]$ of vertices of
level $k$. $A_{[v]}$ is either free or free abelian, and the vertices in $[v]$ do not commute with any vertex that is in the complement of
$[v]$ in $E_k$. 

Recall that in the centerless case,  $A_{E_k}=B*A_{[v]}$, where $A_{[v]}$ is either free or free abelian,
and $B$ is the
group that is generated by the generators that are associated with all the vertices of lower level in $E_k$. With the hyperedge
$E_k$ we add the group of automorphisms of $A_{\Gamma}$ that is constructed as follows.

\roster
\item"{(1)}" We look at generators in the complement of the star of $[v]$, $st([v])$, in $\Gamma$: $\Gamma \setminus st([v])$. 
Let $\Gamma^{[v]}_1,\ldots,\Gamma^{[v]}_m$ be 
the connected components of that complement that are not single vertices. Some of these components contain vertices of lower levels in the hyperedge $E_k$.
Let $S_1,\ldots,S_r$, $1 \leq r \leq m$, be the non-empty subsets of vertices of lower level in $E_k$ that are contained in distinct connected components. Then: 
$A_{E_k}=B*A_{[v]}=B_1*\ldots*B_r*A_{[v]}$, where $B_i$ is the  group that is generated by $S_i$, $1 \leq i \leq r$.

If $A_{[v]}$ is free we first add all the automorphisms of the free product:
$A_{E_k}=B*A_{[v]}=B_1*\ldots*B_r*A_{[v]}$, that conjugate the subgroups $B_i$, $1 \leq i \leq r$ (elementwise). 
The generators of the free factor $A_{[v]}$ are
considered as free generators of the free product.

If $A_{[v]}$ is free abelian we first add all the automorphisms of the free product:
$A_{E_k}=B*A_{[v]}=B_1*\ldots*B_r*A_{[v]}$, that conjugate the subgroups $B_i$, $1 \leq i \leq r$ (elementwise),  i.e., to
elements in the general linear group that is isomorphic to the automorphism group of $A_{[v]}$.

\item"{(2)}" For each connected component $\Gamma^{[v]}_j$, $1 \leq j \leq m$, that contains one of the sets $S_i$, $1 \leq i \leq r$, we conjugate the whole 
group that is associated with a component $\Gamma^{[v]}_j$ by the element that conjugates $B_i$.

Every other component, $\Gamma^{[v]}_j$, that does not contain vertices from the hyperedge $E_k$, and is not a single vertex (i.e., a root), 
is conjugated by an arbitrary element from the group
$A_{E_k}$.

The particular automorphisms that we constructed are automorphisms of the free product:
$$A_{[v]}*A_{\Gamma^{[v]}_1}*\ldots*A_{\Gamma^{[v]}_m}$$
that extend naturally to automorphisms of $A_{\Gamma}$ by defining them to be the identity on the generators that are in $lk([v])$, and on vertices in 
connected components of an isolated (single) vertex in $\Gamma \setminus st([v])$. 
\endroster

Suppose that $E_k$ is a hyperedge with center. Let $[v]$ be the equivalence class of the vertices with highest level in $E_k$. Recall that in the presence of a 
center: $A_{E_k}=Ab \oplus (B*A_{[v]})$, where $Ab$ is in  the center, and $B$ is the (possibly trivial) 
subgroup generated by all the generators in $E_k$ that are in the
complement of $[v]$ and the vertices that commute with the vertices in $[v]$.
In that case we add  automorphisms as follows.

\roster
\item"{(1)}" We start by adding automorphisms similar to the ones that were added for centerless hyperedges. 
We look at generators in the complement of the star of $[v]$, $st([v])$, in $\Gamma$: $\Gamma \setminus st([v])$. 
Let $\Gamma^{[v]}_1,\ldots,\Gamma^{[v]}_m$ be 
the connected components of that complement that are not isolated (single) vertices. 

Some of these  components may contain vertices of lower levels in the hyperedge $E_k$.
Let $S_1,\ldots,S_r$,  be the non-empty subsets of vertices of lower level in $E_k$ that are contained in distinct connected components (if $B$ is trivial
then there are no such subsets of vertices). Then: 
$$A_{E_k}=Ab \oplus (B*A_{[v]})=Ab \oplus (B_1*\ldots*B_m*A_{[v]})$$ 
where $B_i$ is the  group generated by $S_i$, $1 \leq i \leq r$. 

We first add all the automorphism of the free product:
$B_1*\ldots*B_r*A_{[v]}$, that conjugate the subgroups $B_i$, $1 \leq i \leq r$ (elementwise).
If $A_{[v]}$ is free then the automorphisms of the free product regard the generators of $A_{[v]}$ as free generators.
If $A_{[v]}$ is non-cyclic free abelian, then the automorphisms of the free product regard $A_{[v]}$ as a factor, and we add the automorphism group of $A_{[v]}$ which
is 
isomorphic to the general linear group.
We extend such an automorphism to an automorphism of $A_{E_k}$, by defining it to be the identity on the abelian subgroup $Ab$.

\item"{(2)}" As in the case of a centerless hyperedge, for each connected component $\Gamma^{[v]}_j$, $1 \leq j \leq m$, that contains one of the sets $S_i$, 
$1 \leq i \leq r$, we conjugate the whole 
group that is associated with a component $\Gamma^{[v]}_j$ by the element that conjugates $B_i$.

Every other component, $\Gamma^{[v]}_j$, that does not contain vertices from the hyperedge $E_k$,  
is conjugated by an arbitrary element from the group
$A_{E_k}$.

The particular automorphisms that we constructed are automorphisms of the free product:
$$A_{[v]}*A_{\Gamma^{[v]}_1}*\ldots*A_{\Gamma^{[v]}_m}$$
that extend naturally to automorphisms of $A_{\Gamma}$ by defining them to be the identity on the generators that are in $lk([v])$, and on generators that 
are roots that are connected to vertices in $lk([v])$.

\item"{(3)}" We add the remaining transvections of $[v]$. i.e., multiplications of the elements of $[v]$ by  elements from the free abelian group $Ab$, that extend 
to $A_{\Gamma}$ by defining them to be identity on all the generators that are not in $Ab$ nor in $[v]$. 
This is a free abelian group in $Out(A_{\Gamma})$
for each hyperedge with center.
\endroster

We denote the group of automorphisms that is generated by the automorphisms that are associated with all the hyperedges in $\Delta_{\Gamma}$, 
$Aut_1(A_{\Gamma})$.
\endproclaim

\vglue 1.5pc
\proclaim{Theorem 1.3} Let $\Gamma$ be a finite connected graph. Then the image of $Aut_1(A_{\Gamma})$ in $Out(A_{\Gamma})$ is of finite
index in $Out(A_{\Gamma})$.
\endproclaim

\nfp 
Recall that by [La] and [Ser] it suffices to show that the group of automorphisms 
that is generated by $Aut_1(A_{\Gamma})$ and the inner automorphisms of $A_{\Gamma}$ contains all the inversions, 
partial conjugations and transvections of $A_{\Gamma}$.

Inversions of generators are included in the automorphisms of the first type that are associated with the various hyperedges.
Partial conjugations occur when a star of a vertex separates the graph $\Gamma$. They are also contained in the automorphisms of the first type that
are associated with the hyperedges.
   
Transvections occur when there exist two vertices $v,w \in \Gamma$, that satisfy: $lk(v) \subset st(w)$. In that case we can replace 
$v$ by $vw$. If $v$ and $w$ are not connected by an edge, then $lk(v) \subset lk(w)$, so $[v] \leq [w]$. Hence, $v$ and $w$ must be vertices in a
hyperedge in which $v$ is a vertex of highest level. If $[v]=[w]$ the transvections are included in the automorphisms of the free group $A_{[v]}$ that
are included in the automorphisms of the hyperedge that contains $[v]$ as vertices of highest level.
If $[w] > [v]$ then $[v]$ must be an anti-clique and transvections are included with automorphisms of the free product $B*A_{[v]}$ that
is associated with the hyperedge (if it is centerless or not).

Suppose that $v$ and $w$ are connected by an edge in $\Gamma$. Again $[w] \geq [v]$, so $[w]$ is contained in the hyperedge $E_k$ that contains $[v]$ 
as vertices of highest level.
If $[v]=[w]$, then $A_{[v]}$ is a free abelian group and the transvection is included in the general linear group that is associated with $A_{[v]}$ and
with the hyperedge $E_k$.

Suppose that $[w]>[v]$. In that case $[w]$ is a clique,  $E_k$ has center,  and $A_{[w]}$ is contained in $Ab$ that is contained in the center of $E_k$. 
In that case the transvection is contained in the transvections that are associated with the hyperedge $E_k$ that has center.

\line{\hss$\qed$}

If a RAAG is freely decomposable, its flags hypergraph is composed from the flags hypergraphs of its factors, and its 
automorphism group is obtained from the automorphism groups of the factors by adding automorphisms of a free product. 

The flags hypergraph encodes the group of automorphisms of a general RAAG. 
With each hyperedge in the flags hypergraph, and with the 
automorphisms that are associated with the hyperedge, it is possible to associate an action of the RAAG on a simplicial tree, or in case the hyperedge has center,
an action of the RAAG on a simplicial tree in addition to several simplicial actions of the RAAG on lines.

\vglue 1.5pc
\proclaim{Definition 1.4} Let $\Gamma$ be a finite  connected graph, let $A_{\Gamma}$ be its associated  RAAG, and
let $\Delta_{\Gamma}$ be be its flags hypergraph.  
With a hyperedge in $\Delta_{\Gamma}$ we associate an action of $A_{\Gamma}$ on a simplicial tree. If the hyperedge $E_k$ has an abelian associated group,
$A_{E_k}$, we also associate 
with the hyperedge finitely many simplicial actions of the RAAG on lines. 

The associated actions are not faithful in general, and they are all weakly acylindrical (definition 2.1).
We list the possible associated graphs of groups.
\roster
\item"{(1)}" The simplicial actions that are associated with a hyperedge of level 1 can be one of two possibilities.
If the group that is generated by the generators of level 1 is free, 
then we associate with it a graph of groups that contains a bouquet of $s$ loops, where $s$ is
the number of vertices of level 1 in the hyperedge. It contains a vertex for each connected component of $\Gamma \setminus st([u])$ that is not an
isolated (single) vertex, 
where
$[u]$ is the equivalence class of the vertices of level 1 in the hyperedge. 
The stabilizer of such a vertex is the group that is associated with the connected component.
It contains additional $t$ loops that are associated with the connected components in $\Gamma \setminus st([u])$ that are isolated  vertices.

The edge groups in the graph of groups are all trivial.
The kernel of the map from $A_{\Gamma}$ to the fundamental group of the associated graph of groups is normally generated by the generators in $lk([u])$.
The bass-Serre tree that is associated with the graph of groups is the $A_{\Gamma}$-tree that is associated with such a level 1 hyperedge.

Suppose that $[u]$, the equivalence class of level 1, forms a clique with more than one vertex. In that case we start with a  graph of groups that
contains a vertex group
$A_{[u]}$, which is free abelian, and vertex groups that are associated with the 
connected components in $\Gamma \setminus st([u])$, that are not isolated vertices. It contains additional $t$ loops that are associated with the
connected components of $\Gamma \setminus st([u])$ that are isolated vertices. The edge groups are all trivial.
The kernel of the map from $A_{\Gamma}$ to the fundamental group of the associated graph of groups is normally generated by the generators in $lk([u])$.

We further add $s$ simplicial actions of the RAAG on lines, where $s$ is the rank of the free abelian group $A_{[u]}$. The kernel of each of these $s$ actions contains 
all the standard generators of the RAAG, except one of the standard generators of $A_{[u]}$. Hence, the RAAG retracts onto the cyclic subgroup that is generated by the 
standard generator from $A_{[u]}$ that is not in the kernel, and this cyclic subgroup acts simplicially on a line. 
 
\item"{(2)}" The simplicial tree that is associated with a centerless hyperedge of a higher level is similar. Let $[u]$ be the highest level vertices in
the hyperedge. In that case $A_{[u]}$ is free (possibly cyclic). 
The graph of groups that we associate with the centerless hyperedge 
contains a bouquet of $s$ loops, where $s$ is the number of vertices in $[u]$. 
It also contains vertices for each of the connected components of $\Gamma \setminus st([u])$  that are not isolated vertices, where the corresponding vertex 
group is the group that is associated with such a connected component. It contains another $t$ loops, where $t$ is the number of
connected components of $\Gamma \setminus st([u])$ that are isolated vertices. 
All the edge groups are trivial, and the kernel of the map from $A_{\Gamma}$ to
the fundamental group of the graph of groups is identical to the one in part (1). The (simplicial) tree that is associated with such a hyperedge is the 
corresponding Bass-Serre tree.


\item"{(3)}" Suppose that a hyperedge of level more than 1 has center. Let $[u]$ be the highest level vertices in $E_k$.
We look at the vertices in $\Gamma \setminus st([u])$. With this collection of vertices, and their connected components, together with the vertices in $[u]$,
 we associate  a similar graph of groups that we associated
with a centerless hyperedge, depending on whether the vertices of highest level in the hyperedge form a clique or not.

If $A_{[u]}$ is a rank $s$ free group, then with each of the standard $s$ free generators of $A_{[u]}$ there is an associated  loop in the graph of groups,
precisely as in (2). If $A_{[u]}$  is a non-cyclic free abelian group, then the bouquet of $s$ circles in the graph of groups in (2) is replaced by a single
vertex that is stabilized by $A_{[u]}$.

All the edge groups in the constructed graph of groups are trivial. The kernel of the map from $A_{\Gamma}$ to the fundamental group of
the constructed graph of groups is normally generated by the vertices in $lk([u])$. As in (1) and (2), the simplicial tree that is associated with this graph 
of groups is its Bass-Serre tree.
  


If $A_{E_k}$ is abelian, we further add $s$ simplicial actions of the RAAG on lines, where $s$ is the rank of $A_{[u]}$, that are similar to the ones that
were added in case (1) in case $A_{[u]}$ is free abelian.
The kernel of each of these $s$ actions contains 
all the standard generators of the RAAG, except one of the standard generators in $A_{[u]}$. As in (1), the RAAG, $A_{\Gamma}$, 
 retracts onto the cyclic subgroup that is generated by the 
standard generator from $A_{[u]}$ that is not in the kernel, and this cyclic subgroup acts simplicially on a line. 
\endroster 

In all the constructed graphs of groups, the stabilizer of an edge in the fundamental group of the graph of groups (which is a quotient of the RAAG) is trivial. Hence, the
RAAG acts weakly acylinrically (definition 2.1) on the Bass-Serre trees that are associated with the constructed graphs of groups.
\endproclaim

\smallskip
Our next step is to prove that every class of automorphisms in $Out_1(A_{\Gamma})$ 
restricts to outer automorphisms of the fundamental groups 
of the graphs of groups of the first type
that are associated with the various hyperedges that are encoded by these graphs of groups (i.e., the graphs of groups that are not associated with a simplicial action on a line). 
This is important for understanding the structure of the higher rank JSJ
decomposition that is associated with a RAAG, that is obtained from convergent sequences of actions of the RAAG on the trees that we constructed, where the actions
are twisted by automorphisms from $Aut_1(A_{\Gamma})$. 

\vglue 1.5pc
\proclaim{Proposition 1.5} Let $\Gamma$ be a finite connected graph, let $A_{\Gamma}$ be its associated RAAG, and
let $\Delta_{\Gamma}$ be be its flags hypergraph.  With the hyperedges of $\Delta_{\Gamma}$ we have associated graphs of groups of one or two types.

Let $E_k$ be a hyperedge of level $k$ in $\Delta_{\Gamma}$, and let $\tau \in Out_1(A_{\Gamma})$. Suppose that the group that is associated with $E_k$, 
$A_{E_k}=Ab \oplus (B*A_{[u]})$, is not
abelian (i.e., $B$ is non-trivial). Then $\tau$ restricts to an outer automorphism of the 
fundamental group of the graph of groups of the first type that is associated with $E_k$, and its restriction is encoded by this graph of groups. 
\endproclaim 

\nfp 
We start with hyperedges of level 1, $E_1$. Let $[u]$ be the equivalence class of the vertices (of level 1) in $E_1$. Since we assumed that the group that is associated with $E_1$
is not abelian, $A_{[u]}$ is a non-abelian free group. 

By proposition 3.2 in [CV], automorphisms in $Aut_1(A_{\Gamma})$ preserve the conjugacy class of $A_{[u]}$. 
Hence, given 
an automorphism $\tau \in Aut_1(A_{\Gamma})$, we can compose it with 
an inner automorphism, and assume that it preserves the subgroup $<[u]>$. 

$E_1$ is of level 1, hence, $lk([u])$ contains more than a single vertex.
Let $\Theta_{E_1}$ be the graph of groups that is associated with $E_1$ in definition 1.4, and let $G_{E_1}$ be its fundamental group.  
Inversions clearly restrict to  automorphisms 
that are associated with $\Theta_{E_1}$. If $v \in [u]$, then any partial conjugation  in $v$ preserves
the free product in $\Theta_{E_1}$. Partial conjugations by elements
in $lk([u])$ act trivially on $G_{E_1}$, since these generators are mapped to the identity in $G_{E_1}$. 

If $v$ is a vertex
in a connected component of $\Gamma \setminus st([u])$, then a partial conjugation in this vertex may conjugate $<[u]>$ and possibly few other 
factors that are associated with connected components in $\Gamma \setminus st([u])$ by $v^{ \pm 1}$, and leave other such factors unchanged. These are automorphisms 
of free products that are encoded by $\Theta_{E_1}$.

Suppose that $s$ and $v$ are vertices, and a transvection sends $v$ to $sv$. If $v$ is in $E_1$, $s$ must be in $E_1$ as well since vertices in
$E_1$ are maximal. Hence, the transvection is
encoded by $\Theta_{E_1}$.
If $v \in lk[u]$, $s$ must be in $lk([u])$ as well,
so the normal closure of $lk([u])$ does not change, and this transvection does not affect the group $G_{E_1}$, and it doesn't have
an effect on compositions with other automorphisms.

If $v$ is in a connected component of $\Gamma \setminus st([u])$ that is not a single vertex, then either $s$ is in the same component, or $s$ is in
$lk([u])$. Elements in $lk([u])$ are mapped to the identity in $G_{E_1}$. A transvection within a component of $\Gamma \setminus st([u])$ clearly 
extends to an automorphism of a free product that is encoded by $\Theta_{E_1}$. Finally $v$ can be the single vertex in a connected component of
$\Gamma \setminus st([u])$. Such a vertex is represented by a loop in $\Theta_{E_1}$, so any transvection that sends $v$ to $sv$ for some $s$, extends
to an automorphism that is encoded by $\Theta_{E_1}$.

\smallskip
Suppose that for some $k>1$, $E_k$ is a hyperedge for which $A_{E_k}$ is non-abelian. Let $\Theta_{E_k}$ be the graph of groups of the first type that is 
associated with $E_k$ in definition 1.4, and let $G_{E_k}$ be the fundamental group of $\Theta_{E_k}$. The kernel of the map from $A_{\Gamma}$ to
$G_{E_k}$ is normally generated by the generators in $lk([u])$.
 
Let $[u]$ be the equivalence class of the vertices of level $k$ in $E_k$. By the
analysis of the restriction of a general automorphism to the group that is associated with such a hyperedge, the group that is associated with a centerless
hyperedge, $G_{E_k}=Ab \oplus (B*A_{[u]})$, is mapped to a conjugate by every automorphism in $Aut_1(A_{\Gamma})$ (cf. the proof of proposition 3.2 in [CV]).
 
Inversions clearly restrict to  automorphisms 
of $G_{E_k}$ that are encoded by the free product that is associated with $\Theta_{E_k}$. If $v \in E_k$, then any partial conjugation  in $v$ preserves
the free product in $\Theta_{E_k}$, i.e., it may conjugate by $v ^{\pm 1}$ some of the factors and leave others unchanged. Partial conjugations by elements
in $lk([u])$ act trivially on $G_{E_k}$, since these generators are mapped to the identity in $G_{E_k}$. 

If $v$ is a vertex
in a connected component of $\Gamma \setminus st([u])$ that does not contain vertices from $E_k$, then a partial conjugation in this vertex may 
conjugate $<[u]>$ and possibly few other 
factors that are associated with connected components in $\Gamma \setminus st([u])$ by $v^{ \pm 1}$, and leave other such factors unchanged. These are automorphisms 
of free products that are encoded by $\Theta_{E_k}$. 

Suppose that $v$ is a vertex
in a connected component of $\Gamma \setminus st([u])$ that contains vertices from $E_k$ that are not in $[u]$. Then a partial conjugation in this vertex may 
act as an automorphism on the group that is associated with its connected component. It may 
conjugate some of the elements in $[u]$  and possibly few other 
factors that are associated with connected components in $\Gamma \setminus st([u])$ by $v^{ \pm 1}$, and leave other such factors unchanged. 
These are automorphisms 
of free products that are encoded by $\Theta_{E_k}$. 

Suppose that $s$ and $v$ are vertices, and a transvection sends $v$ to $sv$. If $v$ is in one of the connected components of
$\Gamma \setminus st([u])$ that is not a single vertex, then $s$ must be from the same component or from $lk([u])$. Hence, such a transvection
preserves $\Theta_{E_k}$. If $v \in [u]$ then $s$ has to be from $E_k$. Such a transvection is an automorphism that preserves the  free product that
is encoded by $\Theta_{E_k}$.

If $v \in lk[u]$, and $v$ is in the center of $A_{E_k}$, then $s$ must be from the center of $G_{E_k}$ as well, and $v$ and $sv$ are in the kernel of the map from 
$A_{\Gamma}$ onto $G_{E_k}$.
Suppose that $v \in lk([u])$ and $v$ is not in the center of $A_{E_k}$, 
hence, $v \notin E_k$. In that case, $s$ must be from the center of $G_{E_k}$, and again, $v$, $s$ and $sv$ are in the kernel of the map from $A_{\Gamma}$
onto $G_{E_k}$.

Suppose that  $v$ is in a connected component of $\Gamma \setminus st([u])$ that is not an isolated  vertex. Then either $s$ is in the same component, or $s$ is in
$lk([u])$. Elements in $lk([u])$ are mapped to the identity in $G_{E_k}$. A transvection within a component of $\Gamma \setminus st([u])$ clearly 
extends to an automorphism of a free product that is encoded by $\Theta_{E_1}$. 
Finally, $v$ can be the single vertex in a connected component of
$\Gamma \setminus st([u])$. Such a vertex is represented by a loop in $\Theta_{E_k}$, so any transvection that sends $v$ to $sv$ for some $s$, extends
to an automorphism that is encoded by $\Theta_{E_k}$. 

\line{\hss$\qed$}

Proposition 1.5 proves that for hyperedges $E_k$  with non-abelian group, $A_{E_k}$,  outer automorphisms in $Out_1(\Gamma)$ restrict to 
outer automorphisms of $G_{E_k}$ that are encoded by the free product in $\Theta_{E_k}$. Transvections twist the simplicial actions of $A_{\Gamma}$ on
lines that are associated with $E_k$ in definition 1. However, any limit group that is obtained from convergent sequences of twisted actions
of $A_{\Gamma}$ on lines, has to be free abelian, so the family of twisted actions on lines is not difficult to analyze.

If $A_{E_k}$ is abelian, then transvections by elements from $A_{[u]}$ do effect elements from $lk([u])$, so the kernel from 
$A_{\Gamma}$ to $G_{E_k}$ is not preserved by automorphisms from $Aut_1(A_{\Gamma})$. Still, standard generators of $A_{\Gamma}$ that are elliptic in 
$\Theta_{E_k}$ remain elliptic after twisting by automorphisms from $Aut_1(A_{\Gamma})$, so again it is not difficult to analyze the family
of twisted actions of $A_{\Gamma}$ on $\Theta_{E_k}$ also in the case of $A_{E_k}$ abelian. 
 
Maximal decompositions that one obtains from convergent sequences of actions of $A_{\Gamma}$ on the trees that we constructed, where the actions 
are twisted by
automorphisms from $Aut_1(A_{\Gamma})$, are used to construct the higher rank JSJ decomposition of the RAAG. 
This may serve as a motivation or an example for the construction
of a higher rank JSJ decomposition for more general HHGs.

\vglue 1.5pc
\centerline{\bf{\S2. A higher  rank Makanin-Razborov diagram of a product}}
\medskip

In the first  section we associated a flags hypergraph with every RAAG, associated  actions
of the RAAG on one or two types of simplicial trees with each hyperedge in the flags hypergraph, 
and, hence, obtained an action of the RAAG on
products of these (finitely many) simplicial  trees. 


In the next papers in the sequence we generalize these constructions to obtain what we view as a  higher rank JSJ decomposition
of an  HHG that satisfies some further natural assumptions (that hold in the case of the 
mapping class group). In this paper we construct the first step in the construction of the higher rank JSJ decomposition,
what we view as a higher rank Makanin-Razborov diagram.

The Makanin-Razborov diagram was originally constructed (over a free group) to encode the set of solutions to a system of equations
(a $variety$) over a free group, which is equivalent to encoding the set of homomorphisms from a given f.p.\  group (or a f.g.\  group 
if the system of equations is infinite) into a free group (see [Se1]). It was later generalized to study homomorphisms and solutions to systems of equations
over various other families of groups (e.g. [Ja-Se],[Gr-Hu]).

The construction that we present in this paper allows one to construct a higher rank Makanin-Razborov diagram that encodes all the homomorphisms
from a given f.p.\ group into a HHG (that satisfies some mild technical conditions). We believe that the existence of such a diagram (that encodes
homomorphisms) will be used in the near future to generalize properties of hyperbolic groups to HHGs.

However, the whole point of this sequence of papers is to apply techniques and objects that were originally defined and used to study
varieties and more generally first order formulas, to study automorphisms of HHGs. Note that although  homomorphisms from a f.p.\ group can be
identified with a variety, automorphisms are transcadental and in general are not definable (by a first order formula). 

Nevertheless, our approach for studying automorphisms of HHGs is based on applying concepts, objects and techniques that were originally
designed to study first order formulas, to study automorphisms of HHGs. The approach also enables one to construct canonical objects that
encode automorphisms of HHG from 
preliminary constructions that are far from being canonical. We hope that some of these concepts will be generalized to other
families of groups, and to other objects that are not (first order) definable.

\medskip
For presentation purposes we first construct such a higher rank MR diagram in case the HHS  is quasi-isometric to a product of (finitely many)
unbounded hyperbolic spaces, and consider more general HHGs in the next section. Hence, to follow our construction in this section the reader does not need to
be familiar with the definition and the properties of an HHS.

Let $X$ be the product space: $X=V_1 \times V_2 \times \ldots \times V_m$, where each of the spaces $V_j$ are unbounded $\delta$-hyperbolic, for some $\delta>0$.
In the rest of this section we call the spaces $V_j$, $1 \leq j \leq m$, the factors of $X$.

Suppose that a f.g.\ group $G$ acts on $X$ isometrically, properly discontinuously and cocompactly. Let $\pi_j:X \to V_j$, $1 \leq j \leq m$, be the
natural projections. We further assume that $G$ permutes the factors $V_j$, and that 
for every index $j$, $1 \leq j \leq m$, every  $x,y \in X$ and every $g \in G$: 
$$d_{V_j}(\pi_j(x),\pi_j(y))=d_{V_{g(j)}}(\pi_{g(j)}(gx),\pi_{g(j)}(gy)).$$ Hence, a finite index subgroup of $G$ that preserves the factors $V_1,\ldots,V_m$,
acts isometrically on each of the factors $V_j$, $1 \leq j \leq m$.


Note that our assumption on the action of a group $G$ on the product space $X$ coincides with the definition of a group action on an HHS in
the special case of a product space.

In [BHS1] it is proved that a proper co-compact action of a group on a HHS $X$ guarantees that the action of the group on the hyperbolic 
space $CS$ that is associated with the highest complexity space $S$ of $X$ is acylindrical (Corollary 14.4   in [BHS1]). This is an
important property of the action when $CS$ is not a bounded space. 

Acylindricity plays an essential role in our construction of a higher rank JSJ decomposition of an HHG, and in particular in the case of a f.g.\ group that
acts on a product space. In the construction we will need
a weaker form of acylindricity (in the case of a general HHG we will assume this weaker form of acylindricity on all the actions on 
the (unbounded) hyperbolic domains that are associated with the HHS, and not just on the action on the domain that is associated with the highest level one). 
Unfortunately, such weak acylindricity is 
not part of the definition of an 
HHS and it doesn't always hold (e.g., the Burger-Mozes groups [Bu-Mo]). Hence, we need to add it as an additional assumption.

\vglue 1.5pc
\proclaim{Definition  2.1}  Let $Y$ be a metric space and let $G$ act isometrically on $Y$. We say that $G$ acts $weakly$ $acylindrically$ on $Y$ if there exists
$\rho>0$ such that for every $\epsilon >0$, there exist $R,N>0$, so that for every $x,y \in Y$, $d_Y(x,y) >R$, 
there exist at most $N$ elements: $g_1,\ldots,g_k \in G$, 
$k \leq N$, so that if $d_Y(x,gx)<\epsilon$ and $d_Y(y,gy)<\epsilon$, then $g=g_ju$ for some $1 \leq j \le k$, and for every $z \in Y$, $d_Y(z,uz)<\rho$.
\endproclaim

Note that the weak form of acylindricity that we will assume (definition 2.1) is a
coarsification of the WWPD condition of Bestvina, Bromberg and Fujiwara [BBF]. 

\vglue 1.5pc
\proclaim{Definition  2.2}  Let $X$ be a product space $X=V_1 \times \ldots \times V_m$ and let $G$ acts isometrically, properly discontinuously and 
cocompactly as we defined above.

Recall that $G$ permutes the projection spaces $V_1,\ldots,V_m$.  
We say that $G$ acts $strongly$ $acylindrically$ on $X$, if:
\roster
\item"{(1)}"  for every factor $V_j$, $1 \leq j \leq m$, that is not quasi-isometric to a real line,
the (set) stabilizer of $V_j$ in $G$, $stab(V_j)$, acts weakly acylindrically (definition 3.1) on $V_j$.

\item"{(2)}" if $V_j$ is quasi-isometric to a real line, the set stabilizer of $V_j$ modulo the kernel of the action of the set stabilizer on $V_j$, acts
acylindrically on $V_j$.  
\endroster
\endproclaim

In case that a hyperbolic space $V_j$ is not quasi-isometric to a real line, 
an element $u$ that acts quasi-trivially on $V_j$, i.e., $d_{V_j}(z,uz) < \delta$ for every $z \in V_j$,
implies that $z$ is not a loxodromic element. Also, in this case the set of elements that act quasi-trivially on $V_j$ form a normal subgroup of the set stabilizer
of $V_j$.
These properties are important for the analysis of weakly acylindrical actions, and if $V_j$ is quasi-isometric to
a real line, these properties are not always true for weakly acylindrical actions. 

This is why we have to treat the case of spaces that are quasi-isometric to
real lines separately in definition 2.2. However, it is possible to strengthen the weakly acylindrical assumptions in different ways, that still enable one to
analyze the actions of groups on these spaces (note that $[Isom(R):R]=2$ so the algebraic structure of groups that act isometrically on a line are pretty well 
understood).

\smallskip
By our assumptions, $G$ permutes the factors $V_j$, $j=1,\ldots,m$. Hence, there exists a subgroup of finite index, $\hat H$,  that fixes
these hyperbolic spaces. We take $H$ to be the intersection of all the subgroups in $G$ that have the same index as $\hat H$. 
By our assumptions $H$ acts weakly acylindrically on each of the factors $V_j$, 
and $Aut(G)$ acts on $H$.

In order to construct a higher rank JSJ decomposition for $G$, we start by viewing automorphisms of $G$ as homomorphisms, or rather
as quasimorphisms,  and associate a 
(finite) higher rank Makanin-Razborov diagram with $Aut(G)$. Once we have a finite MR diagram, we use the properties of automorphisms, and
in particular the ability to compose them, in order to construct the higher rank JSJ decomposition in the next papers in this sequel.

We look at all the sequences of automorphisms:
$\{\varphi_s \}_{s=1}^{\infty}$ in $Aut(G)$ that are distinct in $Out(H)$.

We fix a generating set $h_1,\ldots,h_{\ell}$ of $H$. Given an automorphism $\varphi_s$, for each index $j$, $1 \leq j \leq m$, there exists 
a point $x^j_s \in V_j$, for which:
$$\max_i \, d_{V_j}(x^j_s,\varphi_s(h_i)(x^j_s))  \ < \ 
1+\inf_{x \in V_j} \max_i \, d_{V_j}(x,\varphi_s(h_i)(x)).$$

Since $X$ is a products of the factors $V_j$, there exists a constant $c>0$, and for each index $s$ a point $x_s \in X$,
such that for each index $j$, $1 \leq j \leq m$, $x_s$ projects to a point that is in a $c$ neighborhood of $x^j_s$ in $V_j$. 

Since $H$ acts cocompactly on $X$, by composing each of the automorphisms in the sequence with an inner automorphism (that depends on the index $s$), 
we may assume that for all $s$: $x_s=x_0$, some  fixed point in $X$.

Since the automorphisms $\{\varphi_s\}$ are distinct in $Out(H)$, there must exist at least one index $j$, for which the sequence:
$\max_i \, d_{V_j}(\pi_{V_j}(x_0),\varphi_s(h_i)(\pi_{V_j}(x_0)))$ is unbounded. Hence, we can extract a subsequence of the automorphisms $\{\varphi_s\}$,
for which these displacements increase to infinity for some $V_j$, $1 \leq j \leq m$. Following [Gr-Hu] we call such a sequence a $divergent$ sequence.

\smallskip
We fixed a sequence of automorphisms in $Aut(G)$ that are distinct in $Out(H)$. We pass to a subsequence for which for some (non-empty set of) indices $j$, 
$1 \leq j  \leq m$,
the projection of the sequence to $V_j$ is divergent, and for the other indices $j$, the sequence is bounded. i.e., for the last indices $j$: 
$\max_i \, d_{V_j}(\pi_{V_j}(x_0),\varphi_s(h_i)(\pi_{V_j}(x_0)))$ is bounded. 

We look at those indices for which the sequence is divergent, and fix such an index $j$.
We can extract a subsequence for which the actions of $H$ on the hyperbolic space $V_j$, twisted by the subsequence of
automorphisms $\{\varphi_s\}$,  converges after rescaling to
a faithful action of a limit group $L$ on a real tree $Y$ (the notion of a convergent subsequence of a divergent sequence is somewhat problematic, 
but we preferred to keep the
existing terminology). 

It is important to note  that limit groups over free groups, that were obtained from sequences of homomorphisms from a given f.g.\ group into a free group, were
defined in two equivalent ways. The first uses Gromov-Hausdorff convergence into actions on real trees, and the stable kernel is then the normal subgroup
that acts trivially on the limit tree. The second definition uses algebraic  convergence, i.e., sequences in which every element of the domain is either 
stably trivial (i.e., mapped to the identity by all except for finitely many homomorphisms), or stably non-trivial. In the second case the stable kernel is
the collection of stably trivial elements.

The two definitions are identical in studying limit groups over free or torsion-free 
hyperbolic groups, but they are distinct in our setting. Since we study automorphisms,
the algebraic definition of a limit group forces all the limit groups to be isomorphic to the original group (automorphisms are injective so the 
algebraic stable kernel is trivial). However, when we use Gromov-Hausdorff convergence to define limit groups (of projections of divergent sequences), 
the stable kernel are those elements that act trivially on the limit (real) trees, and these are not trivial in general. Hence, in order to be able
to apply the limit group
machinery, in the sequel  we will always use the geometric definition of limit groups (Gromov-Hausdorff convergence and the associated stable kernel).

Since the original action of $H$ on $V_j$ was assumed weakly acylindrical, the action of the limit group $L$ on the limit tree $Y$
can be analyzed using the results of [Gu]. 
Applying the work
of [Gr-Hu] (theorems 5.15 and 5.20 in [Gr-Hu]),  since $L$ was constructed from a divergent sequence,  $L$ has a non-trivial virtually abelian JSJ decomposition. 

With the convergent subsequence of automorphisms, for which the sequence of actions converges into the limit group $L$, we want to associate not just
a JSJ decomposition, but also a $resolution$. i.e., 
a finite descending chain of (strict and proper) epimorphisms of limit groups, such that the sequence of homomorphisms that are associated
with the terminal limit group 
in the descending sequence are uniformly bounded, i.e., do not have a divergent subsequence.   

To get such a resolution, and afterwards use these resolutions to construct a Makanin-Razborov diagram, 
we use the techniques that
appear in [Ja-Se] to construct a Makanin-Razborov diagram over free products, together with arguments that are used and appear in [Gr-Hu] for acylindrical actions. 
Note that some of the
main results of [Gr-Hu] are not
applicable in our setup, since they extensively use the equationally Noetherian assumption that we chose to omit.

To construct a resolution, we start with a divergent sequence of automorphisms, $\{\varphi_s\}$, for which the actions of $H$ on $V_j$, twisted by 
$\{\varphi_s\}$, converges into an action of a limit group $L$ on a real tree $Y$. Let $\eta:H \to L$ be the canonical quotient map.
We say that an element $u \in L$ is $stably$ $elliptic$, if for some (hence, any) $h \in H$ for which $\eta(h)=u$, there exists some index $s_h$, such that for all
$s>s_h$, $\varphi_s(h)$ is elliptic (when acting on $V_j$). We say that an element $u \in L$ is $stably$ $non$-$elliptic$ 
if there exists $s_h$ such that for all $s>s_h$, $\varphi_s(h)$ is
loxodromic.

We say that  an element $u \in L$ is $stably$ $bounded$ if for some (hence, any)
element $h \in H$ for which $\eta(h)=u$, the traces of the sequence of elements, $\{\varphi_s(h)\}$, when acting on $V_j$, are bounded. That is:
$\min_{x \in V_j} \, d_{V_j}(\varphi_s(h)(x),x)$ is a bounded sequence. Note that the actual bound may depend on the choice of the element $h$ for which $\eta(h)=u$,
but not the boundedness of the corresponding sequences. Note that every stably elliptic element is stably bounded.

Given a divergent sequence $\{\varphi_s\}$, we apply the diagonal argument and pass to a subsequence for which any element $u \in L$ is either:
\roster
\item"{(1)}"  stably elliptic.

\item"{(2)}"  stably bounded and stably non-elliptic.

\item"{(3)}"   $u \in L$ is not stably bounded, and for every
$h \in H$ for which $\eta(h)=u$, the sequence $\varphi_s(h)$ has no subsequence with bounded traces.
\endroster 

We pass to a subsequence of the original sequence of automorphisms, $\{\varphi_s\}$, that satisfies the above trichotomy. 
We denote the collection of stably bounded elements in $L$ by $B_L$, and the collection of stably elliptic elements by $E_L$ 
(note that these are subsets and not subgroups
in general). 

With the action of $L$ on the limit tree $Y$ we associate a decomposition (graph of groups) of $L$, $\Delta_L$,  
using the analysis of the action as it appears in [Gu]
(see theorem 5.1 in [Gu]). Since we assumed that the action of $H$ on each of the subspaces $V_j$ is weakly acylindrical,
 lemma 4.7 in
[Gr-Hu] proves that the stabilizers of tripods in $Y$ are uniformly finite, the stabilizers of unstable segments are uniformly finite, 
and the stabilizers of non-degenerate segments in $Y$ are uniformly finite by abelian. 

Therefore, every segment in $Y$ is piecewise stable. i.e., every non-degenerate segment $I$ in $Y$ can be divided into finitely many subsegments , $I_j \subset I$,
 such that the stabilizer of $I_j$ is equal to the stabilizer of any non-degenerate subsegment of $I_j$. This enables us to apply the results of [Gu] 
(theorem 4.1 in [Gu]), and deduce that either $L$ splits over the (finite) stabilizer of a tripod, or the action of $L$ on $Y$ decomposes into a graph of actions
(see theorem 4.1 in [Gu] for these notions and conclusion).

If $L$ splits over the stabilizer of a tripod, or if the graph of actions contains an edge with a finite stabilizer, $L$ splits over a uniformly finite group.
In that case we split $L$ over that finite group, and continue with each of the f.g.\ vertex groups. 
Note that all the infinite point stabilizers in $L$ can be conjugated into these f.g.\ vertex groups, and so are in particular the stably elliptic
and stably bounded elements.   

The f.g.\ vertex groups (from the splitting of $L$ along  uniformly finite groups), being subgroups of $L$, act on the
limit tree $Y$. Once again, if one of these f.g.\ groups split over the stabilizer of a tripod in $Y$, or if its associated graph of actions contains
an edge with a (uniformly) finite edge groups, we split the vertex group along the finite subgroup and continue with each of the obtained f.g.\ vertex groups.

From these splittings of f.g.\ vertex groups along uniformly finite subgroups, we get an increasing sequence of graph of groups decompositions of the original
limit group $L$ along uniformly finite subgroups. Such an increasing sequence of graphs of groups has to terminate by a theorem of Linnel [Li], or alternatively
by acylindrical accessibility [We1]. Hence, after finitely many steps of possible refinements, we get a graph of groups decomposition, $\hat \Delta_L$,
 of the original
limit group $L$ with uniformly finite
edge groups, and such that the action of each of its vertex groups on the limit tree $Y$ is either degenerate, or it splits as a graph of (geometric)
actions as it appears in theorem 4.1 in [Gu].
  
By construction each of the vertex groups in $\hat \Delta_L$ acts on the real tree $Y$, and this action splits as a graph of geometric actions according to
theorem 4.1 in [Gu]. With this splitting we can associate a graph of groups decomposition with each f.g.\ vertex group in $\hat \Delta_L$. Since the edge
groups in $\hat \Delta_L$ are all (uniformly) finite, they are contained in vertex groups in these last graphs of groups, so we can use the graph of groups
that are associated with the vertex groups in $\hat \Delta_L$, to further refine it and obtain the graph of groups decomposition $\Delta_L$ of the
limit group $L$.

$\Delta_L$ is a JSJ like decomposition. i.e., it  contains uniformly finite, and uniformly finite by abelian edge groups, and vertex groups that are either:
\roster
\item"{(1)}" rigid (vertex stabilizers that are point stabilizers in the real tree $Y$).

\item"{(2)}"  uniformly finite by abelian, that act axially on the real tree $Y$.

\item"{(3)}" 2-orbifold by uniformly finite (QH), that admit a surface type action on the real tree $Y$.
\endroster

We continue by associating a modular group of automorphisms of $L$ with the decomposition $\Delta_L$, as it appears in section 5.4 of [Gr-Hu]. We replace the 
limit group $L$, that may not be f.p.\ by a sequence of f.p.\ covers $\{M_i\}$ that converge to $L$, and have decompositions (graph of
groups)  with a similar structure
as $\Delta_L$,
their vertex groups are f.p.\  and edge groups are f.p. and uniformly finite by abelian. The QH vertex groups that appear in this sequence of decompositions
are isomorphic to the QH vertex groups in $\Delta_L$. Such a sequence of covers is built and presented in lemma 6.1 in [Re-We] and 
lemma 6.3 in [Gr-Hu].

In the sequel, we say that a map $\nu: G \to Isom(X)$ for some metric space $X$,
 is a quasimorphism if there exists some $q>0$, such that for every $g_1,g_2 \in G$, and every $x \in X$:
$$d_X(x,\nu(g_2^{-1})\nu(g_1^{-1})\nu(g_1g_2)(x)) \ < \ q.$$
Since the limit group $L$ is f.g.\ but may be not f.p.\  the automorphisms $\{\varphi_s\}$ of $H$ from the sequence that converges into
the limit group $L$, may not correspond to  quasimorphisms of the limit (quotient) group $L$. i.e., in general
it is impossible to associate with them  quasimorphisms from $L$ to $Isom(V_j)$. However, with a suffix of the sequence $\{\varphi_s\}$ it is possible
to associate
quasimorphisms of
an increasing sequence of the sequence of f.p.\ covers $\{M_i\}$. i.e., for large index $n$, with the automorphism: $\varphi_s \in Aut(G)$ from
the convergent sequence, there is an associated 
quasimorphism:  $\hat {\varphi_s}:M_{i(s)} \to Isom(V_j)$, and the sequence of indices $\{i(s)\}$ grows to $\infty$. These quasimorphisms can be shortened by  
automorphisms from the modular group of $L$, that lift to automorphisms from the modular groups of the covers $\{M_i\}$.

We denote the shortened quasimorphisms  $\psi_s=\varphi_s \circ f_s$, where $f_s$ is a shortening automorphism from the modular group of $M_{i(s)}$. If the sequence
$\{\psi_s\}$ is a divergent sequence, then a subsequence 
of the quasimorphisms $\{\psi_s\}$ converges to an action of a limit group $L_2$ on a real tree $Y_2$, where $L_2$ is a quotient
(but not necessarily a proper quotient) of the limit group $L$.

Let $\Delta_{L_2}$ be the virtually abelian decomposition that is associated with the action of $L_2$ on $Y_2$, precisely as we associated
the virtually abelian decomposition $\Delta_L$ with the limit group $L$. 
If $L_2$ is a proper quotient of $L$, we continue to the next step, precisely as we did in the construction of the action of $L_2$ on $Y_2$, starting with the
action of $L$ on $Y$. Since the action of $L_2$ on $Y_2$ was obtained from shortening the quasimorphisms $\{\varphi_s\}$ using the modular automorphisms
that are associated with $\Delta_L$, if
$L$ is isomorphic to $L_2$, $\Delta_{L_2}$ contains virtually abelian decompositions that do not appear in $\Delta_L$. i.e., the virtually abelian
decomposition in $\Delta_{L_2}$ give decompositions of rigid vertex groups in $\Delta_L$ that are compatible with $\Delta_L$, or it gives virtually abelian 
decompositions that are hyperbolic w.r.t. $\Delta_L$. 

In both cases it is possible to construct a JSJ type refinement of both $\Delta_L$ and $\Delta_{L_2}$, $\Delta^2_{L_2}$, from which it is possible to
extract both decompositions. Since $\Delta^2_{L_2}$ is a proper refinement of $\Delta_L$, either the decomposition over (uniformly) finite groups
of $L$ has changed, or the modular group that is associated with $\Delta^2_{L_2}$  is strictly bigger than 
the modular group that is associated with $\Delta_L$. Hence, we can now repeat the construction of the action of $L_2$ on $Y_2$, by shortening with
automorphisms from the bigger modular group that is associated with $\Delta^2_{L_2}$.

As long as the shortening quotients are isomorphic to the limit group $L$, we get a sequence of proper refinements of the associated virtually
abelian decompositions. 

\vglue 1.5pc
\proclaim{Lemma 2.3} A sequence of proper refinements of the virtually abelian decompositions of the limit group $L$ terminates
after finitely many steps.
\endproclaim

\nfp By Linnel [Li] or alternatively by acylindrical accessibility [We1], since the decompositions over finite groups along the refinement process, is over 
uniformly finite groups, edges with finite stabilizers can be added finitely many times along the process. 


By lemma 6.2 in [Gr-Hu] a stably non-elliptic virtually abelian subgroup in a non-divergent limit group is virtually cyclic (uniformly finite by cyclic). Hence,
any non-elliptic virtually abelian subgroup in the limit group $L$ is f.g.\ uniformly finite by abelian. 

The iterative process starts with a limit group $L$ and a virtually abelian decomposition $\Delta_L$ and iteratively properly refines it. Since any non-elliptic
virtually abelian vertex group and edge group is f.g.\, the edge groups in the abelian decompositions of $L$ that are constructed along the iterative process,
can be changed (unfolded) only finitely many times along the process.

By proposition 5.14 in [Gr-Hu], the abelian decompositions that are constructed along the process are all $(2,c)$-acylindrical (see [Gu-Le] proposition 7.12). 
Hence, By acylindrical accessibility (see [We1]), there is a bound
on the topological complexity of the orbifolds that are associated with the QH vertex groups in the constructed decompositions, and on the number of edges in these
decompositions. The bound on the combinatorial complexities of the constructed abelian decompositions, together with the finiteness of the numbers  of unfoldings
for each edge, prove that the refinement procedure terminates after finitely many steps.
  
\line{\hss$\qed$}

By lemma 2.3, after finitely many steps, i.e., finitely many proper refinements,  the iterative procedure replaces the limit group $L$ with a proper quotient. 
By theorem 6.8 in [Gr-Hu] every
sequence of proper quotients of limit groups terminates after finitely many steps. Hence, after 
finitely many steps the original sequence, $\{\varphi_s\}$, is replaced by iteratively shortened subsequence that is uniformly bounded, and the limit group
$L$ is replaced by a non-divergent limit group. Altogether, we get a sequence of quotients:
$$ H \, \to \, L \, \to \, L_2 \, \to \, L_3 \, \to \, \ldots \, \to \, L_f$$

Furthermore, the limit groups along the resolution, $L_1,\ldots,L_f$, were constructed as limits of sequences of 
shortened quasimorphisms of the original subgroup $H$.
In particular, the terminal limit group in the resolution, $L_f$, is a limit of a non-divergent sequence of quasimorphisms of $H$. 
Hence, if we denote this non-divergent
sequence of quasimorphisms, $\psi^f_s$, then for any of the fixed set of generators $h_1,\ldots,h_{\ell}$ of $H$, the sequence $d_{V_j}(x_s,\psi^f_s(h_i)(x_s))$
is bounded. Let $bd_f>0$ be a bound on these sequences for all the generators, $h_1,\ldots,h_{\ell}$. To the terminal limit group, $L_f$,
 in the resolution that we constructed from some sequence of automorphisms, $\{\varphi_s\}$, we add this bound $bd_f$, that bounds the displacement 
of the basepoint under the 
images 
of the fixed set of generators of $H$ (images under the terminal sequence of shortened quasimorphisms).

Note that some of the epimorphisms in this sequence are isomorphisms, and some are proper quotients, and the terminal limit group $L_f$ is a 
non-divergent limit group. Also, note that since we gradually enlarged the modular group in case the epimorphisms are isomorphisms, the modular groups that are 
associated with the limit groups  along the resolution are modular groups that are defined using the gradually refined virtually abelian decomposition of the 
corresponding limit group, and not with the virtually abelian decomposition that is associated with the limit group in the resolution.
Hence, the resolution that we obtained is not a strict resolution (see definition 5.9 in [Se1]). In particular, it may be that rigid vertex groups, 
and even edge groups are not mapped
isomorphically by proper epimorphisms. 

\smallskip
Lemma 2.3 starts with a sequence of automorphisms of $H$, for which the sequence of actions of $H$ on the hyperbolic space $V_j$
is divergent for some
$j$, $1 \leq j \leq m$, and constructs from a subsequence of the original sequence a finite resolution, that terminates with a limit group in which
all the elements in the terminal limit group are bounded. We repeat the same construction for all the indices $j$, $1 \leq j \leq m$, for which
the sequence of actions of $H$ on $V_j$ twisted by the automorphisms $\{\varphi_s\}$ is divergent, and for each such $V_j$ we pass to a further
subsequence and construct a resolution. 

As we pointed out, the resolutions that we construct are not necessarily strict (see definition 5.9 in [Se1]), and we will need to work with strict resolutions
to construct the higher rank JSJ decomposition in the next papers. Hence, we modify the resolutions that we constructed.

Let: $$L_1 \, \to \, L_2 \, \to \, \ldots \, \to \, L_f$$ be a resolution, $Res_j$, that is constructed from  a divergent
sequence of automorphisms of the group $H$ (a characteristic finite index subgroup in the HHG $G$), and its associated twisted actions of $H$ on the hyperbolic
space $V_j$. With each limit group $L_i$, $1 \leq i \leq f-1$ there is an associated virtually abelian decomposition, $\Delta_{L_i}$, and a group of modular automorphism
(that is obtained from a gradually refined virtually abelian decomposition).

Starting with the given resolution, we keep the first and the terminal limit groups in the resolution, $L_1$ and $L_f$, 
together with all the limit groups along the resolution for which the quotient
map into them is a proper quotient (and not an isomorphism). 
The virtually abelian decomposition that is associated with each of the limit groups that we kept, except the terminal one, 
is the virtually abelian decomposition that is obtained as a 
common refinement of all the virtually abelian decompositions that are associated with all the limit groups that are isomorphic to the limit group that we kept
along the original resolution. Finally, we obtain a resolution of the form:
$$L_{i(1)=1} \, \to \, L_{i(2)} \, \to \, \ldots \, \to \, L_{i(t-1)} \, \to \, L_{i(t)=f}$$
where all the epimorphisms are proper, except perhaps the terminal epimorphism. We further associate the bound $bd_f$ with the terminal limit group $L_f$
of the new resolution. 

The new resolution is strict by construction (see definition 5.9 in [Se1] for a strict resolution). By section 1 in [Se2], with a strict resolution 
we can (canonically) associate  a completion. 
The completion starts with the terminal limit group $L_f$, as a bottom level, and adds virtually abelian or QH vertex groups according to the abelian decompositions 
that are associated 
with the limit groups along the resolution, from bottom to top. The limit groups $L_{i(1)},\ldots,L_{i(t)}$, where $i(1)=1$ and $i(t)=f$, 
are all naturally and canonically
embedded into the completion, by construction. We denote the completion, $Comp(Res)$.  

Since the terminal limit group $L_f$ may be infinitely presented, so is the completion $Comp(Res)$. In case the completion is not f.p.\
it is not guaranteed that a subsequence of the sequence $\{\varphi_s\}$ extends to a sequence of quasimorphisms of the completion, $Comp(Res)$.
To guarantee the existence of such a subsequence that does extend,
we will need to modify the resolution, or its completion,
and replace it by a f.p.\ $cover$. 

\vglue 1.5pc
\proclaim{Definition 2.4} Let $L_1 \, \to \, L_2 \, \to \, \ldots \, \to \, L_t$ be a strict resolution (with proper epimorphisms, except perhaps the last one)
that is constructed from either a divergent or a bounded
sequence of automorphisms of the group $H$ on one of the spaces $V_j$ (in case the sequence is bounded the resolution has length 1, no quotient maps). With
$L_t$ there is an associated bound $bd_t>0$. 
We say that a resolution:
$M_1 \, \to \, M_2 \, \to \, \ldots \, \to \, M_t$ is a $cover$ of the given resolution if:
\roster
\item"{(1)}" $M_t$ is f.p.\ and there is an epimorphism $\eta: H \to M_1$.

\item"{(2)}" there is an epimorphism $\eta_i:M_i \to L_i$ for every $i$, $1 \leq i \leq t$.

\item"{(3)}" the epimorphisms of the two resolutions together with the epimorphisms $\eta_i$ between their groups form commutative diagrams.

\item"{(4)}" each of the groups $M_i$ is equipped with a virtually abelian decomposition. This decomposition  is a $lifting$ of the virtually abelian decomposition
of the corresponding limit group $L_i$ in the original resolution. i.e., the maps $\eta_i$ map the edge groups 
and the virtually abelian and QH vertex groups in the abelian decompositions of the groups $M_i$ isomorphically onto the edge groups and the virtually abelian and QH
vertex groups in the virtually abelian  decompositions of the limit groups $L_i$. Every rigid vertex group in the virtually 
abelian decompositions of the $M_i$'s
is mapped epimorphically onto a rigid vertex group in the virtually abelian  decomposition of the limit groups $L_i$, by the epimorphisms $\eta_i$.

\item"{(5)}" the cover resolution is a strict resolution with respect to the abelian decompositions that are associated with each of the $M_i$. 

\item"{(6)}" the bound that is associated with the terminal limit group of the cover $M_t$ is the bound that is associated with the terminal limit group
$L_t$, $bd_t$.
\endroster
\endproclaim

Note that since a cover resolution is strict it has a completion. Because the terminal group in a cover is f.p.\ and all the edge groups are f.p.\
(they are virtually f.g.\ abelian), and all the vertex groups that are added along the levels (that either QH or virtually f.p.\ abelian),
the completion of a cover has to be f.p.\ as well.

Following the procedures that appear in the proofs of theorems 24 and 25 in [Ja-Se], we replace the terminal limit group of the resolution, $L_t$,
 with a f.p.\
approximation. The replacement guarantees that the the new completion (obtained by replacing the terminal limit group) is f.p.\ and that   
the QH vertex groups and virtually abelian edge and vertex groups
that appear in the abelian decompositions that are associated with the limit groups along the original resolution remain unchanged in the new cover. 

\vglue 1.5pc
\proclaim{Lemma 2.5} Every resolution  that is obtained from a sequence of convergent automorphisms, $\{\varphi_s\}$, using
the procedure that we presented:
$$L_1 \, \to \, L_2 \, \to \, \ldots \, \to \, L_t$$ 
has a  cover with a f.p.\ completion.
\endproclaim

\nfp Identical to the proofs of theorems 24 and 25 in [Ja-Se].

\line{\hss$\qed$}


A Makanin-Razborov diagram over a free group encodes all the solutions to a system of equations over the free group,
or alternatively all the homomorphisms from a f.g.\ group into the free group. Our goal in this section is to construct a 
(higher rank) Makanin-Razborov diagram, such that every automorphism in $Aut(G)$ factors through it. To achieve that we defined
and constructed cover resolutions.

Replacing a  resolution by a  cover implies that
a tail of the sequence of automorphisms that were used to construct the resolution factor through the (completion of the) cover. 

\vglue 1.5pc
\proclaim{Definition 2.6} Let $M_1 \to \ldots \to M_t$ be a cover of a (strict) resolution $L_1 \to \ldots \to L_t$ that was constructed from a sequence
of automorphisms using the procedure that we presented. Suppose that $bd_t>0$ is the bound that is associated with the terminal limit groups $L_t$ and $M_t$.

We say that an automorphism $\varphi \in Aut(G)$ factors through  the cover resolution
(or alternatively through the cover completion $Comp(Res)$) if $\varphi$ can be written as:
$$\varphi \, = \, \tau \circ \eta_{t-1} \circ \nu_{t-1} \circ \ldots \circ \eta_1  \circ \nu_1$$
where $\nu_i \in Mod(M_i)$,  $\eta_i:M_i \to M_{i+1}$ is the quotient map along the cover resolution, $i=1,\ldots,t-1$, 
and $\tau:M_t \to Isom(X)$ is a quasimorphism, so that
the displacement of the basepoint under $\tau(h_i)$, $i=1,\ldots,\ell$, is bounded by $bd_t$.
\endproclaim

From a sequence of automorphisms $\{\varphi_s\}$ in $Aut(G)$ we passed to a subsequence to construct a resolution: $L_1 \to \ldots \to L_t$. With this resolution
we associated (not canonically and not uniquely) a cover with a f.p.\ terminal group. 
By construction, a tail of the convergent subsequence of automorphisms that were used to construct the given resolution, factors through a given cover
of the resolution.

Note that by our (geometric) definition of a limit group,
saying that an automorphism or a homomorphism factors through a limit group, does not mean that a relation of the limit group
is mapped to the identity by the
homomorphism, but rather that a
relation is mapped by the homomorphism to an element that shifts every point in the target hyperbolic space a distance that is bounded by some constant
(that depends on the hyperbolic space).
Hence, the (completions of the) cover resolutions that we construct do not really encode homomorphisms that factor through them, 
but rather quasimorphisms of the limit
groups (quotients of $H$) that appear along the cover resolution.

Given a divergent sequence, we passed to a convergent subsequence, such that for every index $j$, $1 \leq j \leq m$, we constructed a resolution,
$Res_j$. If for some index $j$, the subsequence is non-divergent, the resolution is a single step resolution. For every index $j$ for which the subsequence
is divergent, the resolution has at least two steps.

We further replace each of the resolutions
$Res_j$, $j=1,\ldots,m$, that were associated with the convergent subsequence and each of the spaces $V_j$, by a cover. 
Hence, the $m$-collection of resolutions that
is associated with a convergent subsequence, and with the spaces $V_1,\ldots,V_m$,  is replaced by an $m$-collection of covers.

Note that we say that an automorphism $\varphi \in Aut(G)$ factors through an $m$-collection of cover resolution, if it factors through each
of the $m$ cover resolutions from the $m$-collection of covers.  
By construction, a tail of the subsequence of automorphisms  that were used to construct the $m$-collection of resolutions, that are associated
with the spaces, $V_1,\ldots,V_m$, factor (in the sense of definition 2.6) through all the $m$-cover resolutions in the $m$-collection of covers.
Hence, a tail of the convergent subsequence of automorphisms factor through the $m$-collection of covers.

In a similar way to  what is proved in [Ja-Se]
for homomorphisms into free products, our next goal is to prove that there exist finitely many $m$-collections of covers, such that every
automorphism in $Aut(G)$ factors through at least one of the finitely many $m$-collections. i.e., that every automorphism in $Aut(G)$
factors through all the $m$
cover resolutions from at least one of the finitely many $m$-collections of covers.

\medskip
We fix a finite generating set of the f.p.\ group $H=<h_1,\ldots,h_{\ell}>$.
We start with all the possible sequences of automorphisms in $Aut(G)$, all their convergent
subsequences, the $m$-collections of resolutions that are constructed from such a convergent subsequence, and the $m$-collections of covers of such
resolutions. 

Note that the terminal group of each cover resolution (in an $m$-collection of covers) is equipped with a bound $bd_t$, on the displacement of the
images of the fixed set of generators of the terminal limit group in the cover under the quasimorphisms that factor through that terminal limit group.
The bounds that are associated with the terminal limit groups of the covers in an $m$-collection  are part of the definition of the automorphisms that 
factor through the $m$-collection (see definition 2.6). Clearly, we may assume that the bounds that are associated with the terminal  groups in 
the $m$-collections of covers are all positive integers.

Since the completion of a cover is f.p.\ and the groups in the resolution are f.g.\ and the virtually abelian decompositions are along f.p.\ edge groups,
the set of (completions of) covers that are decorated by the abelian decompositions and the terminal (positive integers) uniform bounds 
on the displacements of the images of the fixed set of generators is countable.

We order the countable set of $m$-collections of cover resolutions that were constructed from convergent sequences of automorphisms in $Aut(G)$. Note that
each $m$-collection of covers contains the positive integer bounds that are associated with the terminal groups in the $m$-collection of covers.  
If there are automorphisms in $Aut(G)$ that do not factor through the  first $r$ $m$-collections of covers,  we chose
an automorphism $\varphi_r$ that does not factor through the first $r$ $m$-collections of covers.

If no finite subset of $m$-collections of  covers
suffice, we get an infinite sequence of automorphisms, $\{\varphi_r\}$. 
Given this sequence, we can pass to a subsequence that converges into an $m$-collection of resolutions. 

This $m$-collection of resolutions has an $m$-collection of cover resolutions, and this $m$-collection of covers
 must appear in our ordered list of $m$-collections of covers. Suppose that its place in the ordered list of $m$-collection of covers is $r_0$.
Then for large indices $r>r_0$, automorphisms $\varphi_r$ from the convergent subsequence
do factor through this fixed $m$-collection of covers.  But this contradicts the choice of the automorphisms
$\{\varphi_r\}$, since they were supposed not to factor through $m$-collections that appear in the ordered list in the first $r$ places.

\vglue 1.5pc
\proclaim{Theorem 2.7} Let $G$ act properly, cocompactly and strongly acylindrically (definition 2.2) on a product space $X=V_1 \times \ldots \times V_m$. 
Then there exists a finite index characteristic subgroup $H<G$, and  finitely many $m$-collections of cover resolutions, such that every 
automorphism $\varphi \in Aut(G)$ factors through at least one of the finitely many $m$-collections of cover resolutions.

Furthermore, if $Out(G)$ is infinite,
then for at least one index $j$, $1 \leq  j \leq m$, at least one of these finitely many $m$-collections of cover resolutions contains a cover resolution
with at least two steps.
\endproclaim

\smallskip
Theorem 2.7 associates a finite set of $m$-collections of  cover resolutions with the action of $Aut(G)$ on the finite index characteristic subgroup
$H$, in case the space $X$ is a product of finitely many hyperbolic spaces.
We call this finite set of $m$-collections of  cover resolutions, a $higher$ $rank$ $Makanin$-$Razborov$ 
$diagram$ that is associated with the automorphism group $Aut(G)$. Note that we proved its existence using a compactness argument, so
in general the higher rank diagram is not canonical.

\vglue 1.5pc
\centerline{\bf{\S3. A higher  rank Makanin-Razborov diagram of some HHG}}
\medskip

In the first  section we associated a flags hypergraph with every RAAG, associated one or two types of actions
of the RAAG on simplicial trees with each hyperedge in the flags hypergraph, and, hence, obtained  an action of the RAAG
on a product of finitely many simplicial trees. 

In the previous section we studied the automorphism group of a group that acts properly discontinuously and cocompactly on a 
product of hyperbolic spaces, assuming that the action is strongly acylindrical, i.e., that the induced action on each of the factor spaces
is weakly acylindrical (see definitions 2.1 and  2.2). We showed that with such a cocompact action on
a product of hyperbolic spaces, it is possible to associate a (non-canonical) higher rank Makanin-Razborov diagram, that includes  finitely
many $m$-collections of cover resolutions, where $m$ is the number of factors that are associated with the product space. 

Every automorphism of a group that acts properly and cocompactly on the given  space factors through 
at least one of the finitely many $m$-collections of cover resolutions that form the  higher rank MR diagram. Indeed,
the existence of the diagram 
is the first step in the construction of a higher rank JSJ decomposition, and already demonstrate the linkage between automorphisms and their dynamics
and low dimensional topology.

In this section we use the techniques that appear in the previous section to generalize the construction of the higher rank MR diagram
from discrete cocompact actions on products 
of hyperbolic spaces to hierarchically hyperbolic groups (HHG) (see [BHS1], [BHS2], [S] for the definition and basic properties of HHG). 
To be able to apply the techniques that were used in the case of a product, we further assume that the HHG is $colorable$ (or [BBF] colorable).
i.e., that there exists a finite index subgroup
of the HHG, for which the domains of the HHG in an orbit of the finite index subgroup are pairwise transverse. 
This was proved to be true for the mapping class group by
Bestvina-Bromberg-Fujiwara ([BBF],5.8), and is known to hold in quite a few other cases (see e.g. [Ha-Pe] and [DMS]). 

As in the case of product spaces in the previous section, to apply the techniques we further assume some weak acylindricity conditions 
of the actions on the domains of the HHG,
that suffice to guarantee that the actions
of the finite index subgroup on (limit) trees, that are  obtained from associated projection complexes ([BBF], [BBFS]) are weakly acylindrical.  
Again, it is possible to choose 
such sufficient
assumptions  that hold in the case of (HHS that are associated with) mapping class groups, and as we showed in the first section, in the case
of RAAGs. 

\medskip
Let $G$ be an HHG, that acts cocompactly on a HHS $X$. Suppose that $G$ is colorable, i.e., that there exists a finite index subgroup $H<G$,
such that the domains of the HHS $X$ in each of the (finitely many) orbits under the action of $H$, are pairwise transverse (which means that they
are not orthogonal and not nested). 
Let $m$ be the number of 
orbits of domains in $X$ under the action of $H$ (the number of orbits of domains is assumed to be finite for an HHG).

Let $V_1,\ldots,V_m$ be representatives from the distinct orbits of domains in $X$ under the action of $H$. Then pairs of distinct
domains in each  orbit: $H(V_j)=\{hV_j \, | \, h \in H\}$ are transverse. From each transverse collection, $H(V_j)$, we construct a projection complex and
a quasi-tree of
hyperbolic metric spaces using the constructions of Bestvina-Bromberg-Fujiwara [BBF] (see also [BBFS]). 

[BBF] starts with a collection of metric spaces that they denote $\pmb{Y}$, and with a constant $\theta>0$, 
such that for every $Y \in \pmb{Y}$ there is a function:
$$d^{\pi}_Y: (\pmb{Y} \setminus Y) \, \times \, (\pmb{Y} \setminus Y) \ \to [0,\infty)$$
that has the following properties for all the metric spaces $X,Y,Z,W \in \pmb{Y}$:
\roster
\item"{(PC0)}" $d^{\pi}_Y(Z,Z) \, < \, \theta$.

\item"{(PC1)}" $d^{\pi}_Y(X,Z) \, = \, d^{\pi}_Y(Z,X)$.

\item"{(PC2)}" $d^{\pi}_Y(X,Z)+d^{\pi}_Y(Z,W) \, \geq \, d^{\pi}_Y(X,W)$.

\item"{(PC3)}" $\min(d^{\pi}_Y(X,Z),d^{\pi}_Z(X,Y)) \, \leq \, \theta$.

\item"{(PC4)}" for all $X,Z \in \pmb{Y}$, $\# \{Y \, \vert \, d^{\pi}_Y(X,Z)>\theta \, \}$ is finite. 
\endroster
(see section 3.1 in [BBF]). We assumed that each distinct pair of metric spaces from the collection $H(V_j)$ is transverse. We refer to
definition 1.1 of hierarchically hyperbolic spaces (HHS) in [BHS2]. Given two transverse domains $X,Y$ of an HHS $X$,
there is a  projection set $\rho^X_Y \subset Y$ of uniformly bounded diameter $\kappa_0$, and we can assume $\theta>\kappa_0$.

On the set $H(V_j)$ of transverse domains there are naturally defined pseudo metrics:
$$d^{\pi}_Y(X,Z)=diam_Y(\rho^X_Y \, \cup \, \rho^Z_Y).$$
The definition immediately implies (PC0)-(PC2). (PC3) is a special case of ([BHS2],1.8), and (PC4) 
follows from the distance formula for HHS ([BHS2],4.5), which is a generalization of
Masur-Minsky distance formula for the mapping class groups.   

Theorem 4.1  in [BBFS] proves that by replacing the constant $\theta$ with $11\theta$ it is possible to replace the pseudo-metrics
$d^{\pi}_Y$ with pseudo-metrics $d_Y$ that satisfy the axioms (SP1)-(SP5) of [BBFS]. These include the axiom:

\noindent
$(SP3) \ \ $ if $d_Y(X,Z)>\theta$ then $d_Z(X,W)=d_Z(Y,W)$ for all $W \in \pmb{Y} \setminus \{Z\}$.

\noindent
The properties (SP1)-(SP5) will allow us to work with the tools and the conclusions of [BBFS], and in particular with  their
$standard$ $paths$ in projection complexes (lemma 3.1 in [BBFS]). Standard paths are helpful in analyzing projection complexes and
in particular in proving acylindrical properties of a group action on these complexes (e.g., theorems 3.9 and 6.4 in [BBFS]).

After modifying the pseudometric to satisfy properties (SP1)-(SP5) in [BBFS], let $K>3\theta$. With an orbit $H(V_j)$ we associate the projection 
complex $P^j_K$ (see section  3 in [BBF] or [BBFS]). The following guarantees weak acylindricity of the action of $H$ on the projection complex $P_K^j$.

\smallskip
\vglue 1.5pc
\proclaim{Proposition 3.1 (cf. ([BBFS],3.9))} Let $K>3\theta$. Suppose that there exist some positive integers $r,b$, 
 such that  the stabilizer of any standard path of length  $r$ in the projection complex $P_K^j$
has a finite index subgroup of index at most $b$, for which every element in the finite index subgroup
 stabilizes (setwise)  each of the domains in the orbit $H(V_j)$. Then $H$ acts on the projection complex $P_K^j$ weakly acylindrically (definition 3.1).   
\endproclaim

\nfp Identical to the proof of theorem 3.9 in [BBFS].

\line{\hss$\qed$}
 
The weakly acylindrical actions of the finite index subgroup $H$ on the projection complexes $P_K^j$, $1 \leq j \leq m$, enable us to apply the techniques
that were used in the case of a product in the previous section. Hence, with these $m$ weakly acylindrical actions we can associate a higher rank
Makanin-Razborov diagram that encodes the automorphisms in $Aut(G)$.

\vglue 1.5pc
\proclaim{Theorem 3.2 (cf. Theorem 2.7)} Let $G$ act properly and  cocompactly on an HHS $X$ and
suppose that $X$ is colorable. i.e., there exists a finite index subgroup $H<G$, such that there are only  
finitely many orbits of domains of $X$ under the action of $H$ (i.e., $H$ and $G$ are HHG), and   the domains in each orbit of $H$,
$HV_j$, are pairwise transverse. 
wlog.\ we may assume that $H$ is a characteristic finite index subgroup of $G$.
Let $m$ be the number of orbits of domains under the action of $H$. 

For $j$, $1 \leq j \leq m$, let $P_K^j$ be the projection complex that is constructed from the orbit, $HV_j$, using the construction in 
[BBF] and [BBFS]. 
Suppose further that there exist some positive integers $r,b$, 
 such that  the stabilizer of any standard path (see [BBFS] for this notion) of length  $r$ in the projection complexes $P_K^j$, $1 \leq j \leq m$,
has a finite index subgroup of index at most $b$, for which every element in the finite index subgroup
 stabilizes (setwise)  each of the domains in the orbit $H(V_j)$.    

Then there exist finitely many $m$-collections of  cover resolutions, such that the action of $H$
on the projection complexes $P_K^j$, $1 \leq j \leq m$, twisted by any automorphisms in $Aut(G)$, factors through at least one of the $m$-collections of  cover 
resolutions (see definition 2.6 for a twisted action that factors through an $m$-collection). 

Every  cover resolution terminates with a f.p.\  group, with some fixed finite generating set
and a positive integer. The displacement of the basepoint in $P_K^j$, under the image of each of the fixed set of generators under 
all the quasimorphisms from the terminal limit group to the isometry group of $P_K^j$, that are associated with automorphisms 
that factor through the cover resolution, 
are bounded by this positive integer.
\endproclaim

\nfp Unlike the proof of theorem 2.7, to construct the limit trees from (twisted) action of $H$ on the various projection
complexes, $P_K^j$, we choose basepoints that are displaced minimally in each of the projection complexes,
$P_K^j$, separately. i.e., the basepoints in the projection complexes are not necessarily an image of some global basepoint (that is displaced minimally) in the 
ambient HHS $X$. 

Apart from the change in the choice of the basepoints, the proof is identical to the proof of theorem 2.7.

\line{\hss$\qed$}

\smallskip
The proof of theorem 3.2 does not really require properness nor cocompactness of the action of $G$ on the HHS $X$ (but we preferred to keep it as part of our
general theme).
Theorem 3.2 associates a finite set of $m$-collections of  cover resolutions with the action of $Aut(G)$ on $H$, that are constructed from the 
twisted actions of $H$ on the projection complexes, $P_K^j$,  that are associated with the orbits of domains under the action
of $H$, $HV_j$. 

The fixed set of generators of the terminal limit groups in these cover resolutions have uniformly bounded displacements when acting on the projection complexes,
$P_K^j$,
through the quasimorphisms from the terminal limit groups. However,
these bounds are not sufficient to bound the action on the HHS space $X$. In particular, all the cover resolutions in the finite set of  $m$-collections 
can have single levels and still
$Out(G)$ may be infinite.

\medskip
To extract further information on $Aut(G)$ we can continue with two different objects and actions. The two different objects require different assumptions to enable 
us to analyze them using the techniques that appeared in the previous section (the proof of theorem 2.7). 

The first direction is to
start with the actions of $H$ on the  quasi-trees of metric spaces that were constructed in [BBF],
instead of starting with the actions on the projection complexes.  As in [BBF], under our assumptions, the characteristic finite index subgroup  $H<G$ 
is quasi-isometrically embedded in a product of the quasi-trees of metric spaces, and if we assume that the actions of $H$ on the quasi-trees of
metric spaces are weakly acylindrical, we can apply the techniques of the previous section and obtain similar results.

The second direction requires weaker (weakly) acylindrical assumptions, and the statement of its conclusion is more technical. 
In this (second) direction we continue with the 
terminal limit groups in the $m$-collections of cover
resolutions that were constructed from the projection complexes in theorem 3.2, and analyze their actions on 
the product of quasi-trees of metric spaces of [BBF] and [BBFS]. Here we only require that the set stabilizers of the domains act weakly acylindrically on the domains that
they stabilize, and the global action of $H$ on the quasi-trees on metric spaces is not assumed to be weakly acylindrical. In this direction we still construct a higher
rank MR diagram, but instead of containing finitely many $m$-collections of resolutions, the diagram contains finitely many $m$-collections of $hybrid$ 
$resolutions$. Hybrid resolutions
played a central role in the solution of Tarski's problem, in particular in the construction of the $anvil$ and the $developing$ $resolution$ (section 4 in [Se4]).

\vglue 1.5pc
\proclaim{Theorem 3.3 (cf. Theorems 2.7 and 3.2)} Let $G$ act properly and  cocompactly on a HHS $X$. 
Suppose that $G$ is colorable (as in theorem 3.2). i.e., that there exists a finite index characteristic subgroup $H<G$ for which there are 
finitely many orbits of domains of $X$ under the action of $H$, and  the domains in an orbit of $H$ are pairwise transverse.
Let $m$ be the number of orbits of domains under the action of $H$. 

Applying the constructions in [BBF] and [BBFS], with the action of $H$ on $X$ it is possible to construct $m$ actions of $H$ on  
quasi-trees of hyperbolic metric spaces (the domains of $X$ in
an orbit of $H$), $C_K^j$, $1 \leq j \leq m$.
Suppose that the  actions of $H$ on these $m$ quasi-trees of metric spaces are weakly acylindrical (definition 2.1). For those indices $j$ for which
$C_K^j$ is a quasiline, we further assume  that $H$ modulo the pointwise stabilizer of $C_K^j$ acts acylindrically on $C_K^j$ (part (2) in definition 2.2).

Then there exist finitely many $m$-collections of cover resolutions, such that the actions of $H$
on the $m$ quasi-trees of metric spaces, $C_K^j$,  twisted by any automorphism in $Aut(G)$, factor through at least one of the $m$-collections of cover resolutions
(see definition 2.6 for a twisted  action that factors an $m$-collection). 
The quasimorphisms of the terminal groups in the cover resolutions in one of the finite set of $m$-collections, that are associated with
the subset of automorphisms from $Aut(G)$ that factor through that  particular  $m$-collection, are uniformly bounded. i.e., the displacements in $C_K^j$ of the basepoints under
the images of a fixed set of generators are uniformly bounded.

Finally, if $Out(G)$ is infinite, then at least one of the $m$-collections contains a cover resolution with at least two levels.
\endproclaim

\nfp As in the proof of theorem 3.2, the basepoints are chosen separately in each of the quasi-trees of metric spaces, $C_K^j$, to be points with minimal
displacements. The construction of the finitely many $m$-collections that satisfy the conclusions of the theorem is identical to the proof of theorem 2.7
(and theorem 3.2).

Suppose that $Out(G)$ is infinite.
Let $h_1,\ldots,h_{\ell}$ be a fixed generated set of $H$, and
let $\{\varphi_s\}$ be a sequence of automorphism from $Aut(G)$, that belong to distinct classes from $Out(G)$. For each automorphism $\varphi_s$, there is
a point in $X$ that is displaced minimally (up to a global constant) by the elements, $\varphi_s(h_1),\ldots,\varphi_s(h_{\ell})$. Since the action of $H$ on
$X$ is cocompact, after possibly composing the automorphisms $\{\varphi_s\}$ with inner automorphisms, we can  assume that the point that is displaced minimally
(up to a uniform constant) by the images of the generators under the automorphisms $\{\varphi_s\}$ is some fixed point $x_0 \in X$. Since the automorphisms
belong to distinct classes in $Out(G)$, and the action of $H$ on $X$ is proper, the displacements of $x_0$ under the actions, twisted by the automorphisms 
$\{\varphi_s\}$,  of the images of the generators
$h_1,\ldots,h_{\ell}$ is unbounded. After passing to a subsequence of the automorphisms $\{\varphi_s\}$, we may assume that the displacements grow to $\infty$. 

To complete the proof of theorem, assume that all the resolutions in the constructed finite set of $m$-collections of cover resolutions have only a single level.
In that case, there is some bound $b>0$, such that
for each index $s$, and every $j$, $1 \leq j \leq m$, there exists a point $y_s^j \in C_K^j$, such that the displacements of the points $y^j_s$ in $C_K^j$  under the
images of $h_1,\ldots,h_{\ell}$ twisted by the automorphisms  $\varphi_s$ are bounded by $b$. 

We aim to get a contradiction, and we will do that by finding points in $C_K^j$ that are displaced only a bounded distance by the actions of the images of
$h_1,\ldots,h_{\ell}$ twisted by the automorphisms $\{\varphi_s\}$. This will contradict the fact that $x_0$ is displaced minimally (up to a constant), 
and its displacements grow to $\infty$ 
with $s$. 

We fix $j$, $1 \leq j \leq m$. To save notation, we will assume that the identity of $H$ is one of the fixed set of generators $h_1,\ldots,h_{\ell}$.
By the distance formula (theorem 4.5 in [BHS2]), and the construction of the quasi-tree of metric spaces [BBF], $C_K^j$,
the interval $[\varphi_s(h_{i_1})(x_0),\varphi_s(h_{i_2})(x_0)]$, $1 \leq i_1 < i_2 \leq \ell$, is supported on finitely many domains in $C_K^j$. If 
the support is not empty, and since
$C_K^j$ is $\delta$-hyperbolic for some $\delta>0$, the collection of these paths in $C_K^j$ is in some $\delta'>0$-neighborhood of a finite tree,
$T_s^j$, in
$C_K^j$ (where $\delta$ and $\delta'$ do not depend on $s$ nor $j$).

If for some index $j$, $1 \leq j \leq m$, there exists a subsequence of indices $s$ for which the diameters of the trees $T_s^j$ are (globally) bounded,
we pass to this subsequence. Hence, after passing to a subsequence, and changing the order of the indices $j$, we may assume that for $j$, $1 \leq j \leq m'$,
the diameters of the trees $T_s^j$ grow to $\infty$, and for $m'< j \leq m$, the diameters of the trees $T_s^j$ are globally bounded or these trees are empty. 
Note, that since 
the displacements of $x_0$ in $C_K^j$ grow to $\infty$, at least for a single index $j$, the diameters of the trees $T_s^j$ can not be globally bounded.
Hence, $m' \geq 1$.

Since $H$ acts isometrically on $C_K^j$, and for each  $j$, $1 \leq j \leq m'$, and each index $s$,  there exists a point $y_s^j \in C_K^j$ 
that is displaced not more than a distance $b$
by the images of $h_1,\ldots,h_{\ell}$ under the automorphism $\varphi_s$, there is a point $t_s^j \in T_s^j$, that is displaced not more than $b'=b+4 \delta$ by
the images (under $\{\varphi_s\}$) of each of the elements $\varphi_s(h_i)$, $1 \leq i \leq \ell$, in $C_K^j$.    

The points $t_s^1$ are approximately the midpoints  in  intervals, which are the support in $C_K^1$ of the intervals $[x_0,\varphi_s(h_{i_0})(x_0)]$ (the last intervals are in $X$),
for some $i_0$, $1 \leq i_0 \leq \ell$.
By theorem 4.4 in [BHS2] there exist hierarchy paths between any two points in the HHS $X$. Hierarchy paths project to unparameterized quasi-geodesics in each of the domains of
$X$, and are contained in a uniform neighborhood of the hull of the two points in $X$. 

For each $s$, let $\gamma_s$ be a hierarchy path from $x_0$ to $\varphi_s(h_{i_0})(x_0)$. The support of $\gamma_s$ in $C_K^j$ lies within some uniform neighborhood of the 
support of the interval $[x_0,\varphi_s(h_{i_0}(x_0)]$ in $C_K^j$, for $j$, $1 \leq j \leq m$. In particular, for $j$, $m' <j \leq m$, the support of $\gamma_s$ 
is uniformly bounded in $C_K^j$.

Furthermore, there exists a point $x_s^1 \in \gamma_s$, such that the  endpoint of the support in $C_K^1$ of the part of $\gamma_s$ from $x_0$ to $x_s^1$ is
uniformly close to $t_s^1$. i.e., the distance between  the endpoint of the support to $t_s^1$ in $C_K^1$  is globally bounded (for all indices $s$).
Hence, for every index $s$, and every 
$i$, $1 \leq i \leq \ell$, the support of $[x_s^1,\varphi_s(h_i)(x_s^1)]$ in $C_K^1$ is uniformly bounded. Hence, if we replace $x_0$ by $x_s^1$, 
the union of the supports of the intervals:
$[x_s^1,\varphi_s(h_i)(x_s^1)]$, $1 \leq i \leq \ell$, in the spaces $C_K^j$ are unbounded only for a subset of the indices $j$, $2 \leq j \leq m'$. 

Therefore, the union of the supports of the intervals $[x_s^1,\varphi_s(h_i)(x_s^1)]$ are unbounded for a proper subset of the spaces $C_K^j$, for which the union of the
supports of the intervals $[x_0,\varphi_s(h_i)(x_0)]$ were unbounded. Repeating the construction of the points $\{x_s^1\}$ iteratively, we finally get points
$\{x_s\}$ in $X$, such that for every index $s$,  every $j$, $1 \leq j \leq m$, and  every $i$, $1 \leq i \leq \ell$, the union of the supports of the intervals
$[x_s,\varphi_s(h_i)(x_s)]$ in $C_K^j$ are uniformly bounded. Hence, by the distance formula, 
 the distances $d_X(x_s,\varphi_s(h_i)(x_s)]$ are uniformly bounded.

The points $\{x_s\}$ have uniformly bounded displacements under the twisted actions of the elements $h_1,\ldots,h_{\ell}$, and that contradicts the choice of the point 
$x_0$ that was assumed to have minimal displacements (up to some global constant) under these twisted actions, and these displacements grow to $\infty$ 
with $s$ for the sequence of 
automorphisms $\{\varphi_s\}$.

Finally, the contradiction implies that if $Out(G)$ is infinite, then   the finite set of $m$-collections of cover
resolutions that we constructed from $Aut(G)$ must contain a cover resolution with at least two levels.


\line{\hss$\qed$}

Theorem 3.3 assumes that the actions of the finite index subgroup $H<G$ on the $m$ quasi-trees of metric spaces, 
that were constructed from the action of $H$ on $X$, are weakly
acylindrical. It is possible to somewhat relax the conditions on the action of $H$ on $X$, and associate a different type of a higher rank diagram with the action.

The assumption that the action of $H$ on the $m$ quasi-trees of metric spaces is weakly acylindrical, is replaced by the assumption that the action of $H$ on the
$m$ projection complexes is weakly acylindrical and so are the actions of the (set) stabilizers of the domains on the domains that they stabilize. 
The finite $m$-collections
in the higher rank diagram that we construct, are not $m$-collections of cover resolutions, but rather $m$-collections of hybrid resolutions.  

\vglue 1.5pc
\proclaim{Theorem 3.4 (cf. Theorem 3.3)} Under the assumptions of theorem 3.2, suppose that $HV_j$, $1 \leq j \leq m$, are the orbits of the domains in the
HHS $X$ under the action of $H$, and   the set stabilizer in $H$ of each domain $V_j$,
$1 \leq j \leq m$, acts weakly acylindrically on $V_j$ (definition 2.1).  
For those indices $j$ for which the domain $V_j$
is a quasiline, we further assume  that the set stabilizer of $V_j$  modulo the pointwise stabilizer of $V_j$ acts acylindrically on 
$V_j$ (part (2) in definition 2.2).

With the actions of $H$ on $X$ twisted by automorphisms from $Aut(G)$ it is
possible to associate a finite set of $m$-collections of $hybrid$ (cover) $resolutions$.
Each hybrid resolution consists of a pair of resolutions.    
The first  resolution in each pair, is a resolution of the
ambient group $H$, and it is constructed from actions of $H$ twisted by automorphisms from $Aut(G)$, 
on the associated projection complex $P_K^j$, $1 \leq j \leq m$. 

The second part in each pair is a resolution of the terminal limit group of the first resolution. The second resolution
 is composed from finitely many resolutions of
subgroups of the terminal group of the first resolution, which are the intersections of the terminal group with (set) stabilizers of domains in $X$.

Every automorphism in $Aut(G)$ factors through at least one of the finitely many $m$-collections of hybrid resolutions. Note that an automorphism factors through  
an $m$-collection of hybrid resolutions if the action of $H$ on the HHS space $X$ factors through all the $m$ hybrid resolutions in the $m$-collection according
to definition 2.6. 

If $Out(G)$ is infinite, then at least one of the $m$-collections  pair resolutions contains a pair in which at least one of the   resolutions from the pair has
 at least two levels.
\endproclaim

\nfp 
We start with a sequence of automorphisms, $\{\varphi_s\}$, in $Aut(G)$.
Following [BBF] and [BBFS],  from the actions of $H$ on the pairwise transverse orbits of the domains, $HV_j$, $1 \leq j\leq m$, 
it is possible to construct $m$ projection complexes,
$P_k^j$, $1 \leq j \leq m$. 

To construct the top resolution in each hybrid resolution, we repeat what we did in theorem 3.2. We start with the actions of $H$ on the $m$ projection complexes, 
$P_K^j$, $1 \leq j \leq m$, twisted
by the sequence of automorphisms, $\{\varphi_s\}$. The basepoint for each such twisted action is taken to be a point in $P_K^j$ with minimal displacement
under the images of a fixed finite set of generators of $H$. Note that these basepoints are not necessarily the image of a point in the HHS $X$, but
they are chosen separately in each projection complex $P_K^j$. 

By our assumptions, the actions of $H$ on the complexes $P_k^j$ are weakly
acylindrical. Hence, we can apply the construction that was used in the proof of theorem 2.7, pass to a subsequence of the automorphisms
(that are still denoted)
$\{\varphi_s\}$, and construct $m$ resolutions that terminate in limit groups, such that with each terminal limit group there is an
associated graph of groups. Each edge group in this graph of groups has a finite stabilizer.
Each vertex group in the terminal graph of groups is therefore f.g., and there are uniform bounds on the displacement of the basepoint under
the  images of 
a fixed finite set of generators of each vertex group  under the sequence of quasimorphisms that are associated with each terminal limit group.

Let $\{\psi^j_s\}$, $1 \leq j \leq m$, be the sequences of quasimorphisms that converge into the terminal limit groups $L^j_t$, $1 \leq j \leq m$,
in the $m$ resolutions that were constructed from the twisted actions on the projection complexes $P_K^j$, $1 \leq j \leq m$.     

By [BBF], from the actions of $H$ on the orbits of the pairwise transverse domains, $HV_j$, $1 \leq j \leq m$,  it is also possible to construct $m$
quasi-trees of metric spaces that we denote $C_K^j$, in which the unbounded hyperbolic domains $V_j$ are embedded. The sequences
of quasimorphisms $\{\psi^j_s\}$ were obtained from twisted actions of $H$ on the projection complexes $P_K^j$, but they can be viewed as quasimorphisms
into the isometry groups of the quasi-trees of metric spaces $C_K^j$ (since $H$ acts on $C_K^j$, and unless $P_K^j$ or $V_j$ are bounded,
the coarse pointwise stabilizer of the action of $H$ on $P_K^j$ is  the coarse pointwise stabilizer of the action of $H$ on $C_K^j$).
     
We continue in parallel with each of the (finitely many) 
f.g.\ vertex groups in the terminal graph of groups decomposition that is associated with a terminal limit group $L^j_t$,
and with a f.g.\ preimage of such vertex group in $H$. To save notation we will continue to denote such a vertex group $L^j_t$, and the f.g.\ preimage $H$.
 
The uniform bounds on the displacements of the basepoints in $P_K^j$ under the images  
of fixed finite sets of generators  of the limit groups $L^j_t$, 
twisted by the sequence of quasimorphisms, $\{\psi^j_s\}$,  
guarantee that if we choose a basepoint in $C_K^j$ to be in the domain in which the basepoint in $P_K^j$ is located,
then for each of the associated quasimorphisms into the isometry groups of the  quasi-trees of metric spaces, $C_K^j$,
the paths that connect the basepoint to its image by the fixed set of generators of $H$ (twisted by the  quasimorphisms), 
are supported by boundedly many domains
in the quasi-trees $C_K^j$.

If for some indices, $s$ and  $j$, $1 \leq j \leq m$, the union of the supports of the intervals $[x_0, \psi^j_s(h_i)(x_0)]$,
$1 \leq i \leq \ell$,  in $C_K^j$ is empty, we
choose the basepoint, $c_s^j \in C_K^j$ to be a point that is displaced minimally in $C_K^j$ 
(up to some global constant) by the elements $\varphi_s(h_i)$, $1 \leq i \leq \ell$. 

Suppose that for a pair of indices, $s,j$,  the union of the supports of the intervals:
$[x_0, \psi^j_s(h_i)(x_0)]$,
$1 \leq i \leq \ell$,  in $C_K^j$ is not empty. In that case we choose the basepoint, $c_s^j \in C_K^j$,  
to be a point in the union of the supports of these intervals, 
that is displaced minimally (up to some global constant) by the elements: $\psi^j_s(h_i)$, $1 \leq i \leq \ell$, among all the points in the union of the supports
of these intervals in $C_K^j$. 

Note that because $c_s^j$ is a point that is displaced minimally (up to a global constant) along the points in the
supports of the intervals:
$[x_0, \psi^j_s(h_i)(x_0)]$, the difference between its displacement and the displacement of the infimum of the displacements of points in $C_K^j$ by
$\psi^j_s(h_i)$, $1 \leq i \leq \ell$, is globally bounded. Also, since the points $c^j_s$ are contained in the union of the supports of the intervals:
$[x_0, \psi^j_s(h_i)(x_0)]$, in $C_K^j$, and the union of the supports of these intervals in $P_K^j$ are uniformly bounded,
the union of the supports of the intervals:
$[c^j_s, \psi^j_s(h_i)(c^j_s)]$, in $P_K^j$, are uniformly bounded as well.

If for some index $j$, 
the displacements (of the chosen basepoints $\{c_s^j\}$) under the elements $\psi^j_s(h_i)$, $1 \leq i \leq \ell$, 
are not bounded, we can pass to a subsequence that converges into an action of the terminal vertex group $L^j_t$ on a real tree $Y_j$. 

We didn't assume that the actions of $H$ on the quasi-trees of metric spaces, $C_K^j$, are weakly acylindrical, but only that the actions of the domains
(set) stabilizers on the domains that they stabilize are weakly acylindrical. Hence, to construct resolutions from the actions of the quasimorphisms
of $L^j_t$ into the isometry groups of $C_K^j$, and a Makanin-Razborov diagram, we need to analyze not the action of $L^j_t$ on the tree $Y_j$, but rather the
actions of the finitely many set stabilizers  in $L^j_t$ of domains in the orbits of the domain $V_j$, on these domains. 

Note that these set stabilizers of domains need not be f.g.\ but, as we show in the sequel, 
they are f.g.\ relative to finitely many stably bounded subgroups. To start the analysis of these set
stabilizers, and prove their relative finite generation, we start with the construction of the following graph of groups.
 

\vglue 1.5pc
\proclaim{Proposition 3.5} It is possible to pass to a further subsequence of quasimorphisms of $L^j_t$, that we still denote $\{\psi^j_s\}$,
 into the isometry group of the quasitrees of metric
spaces, $C_K^j$, and associate a finite bipartite graph of groups decomposition $\Theta_j$ with the limit group $L^j_t$, for every $j$, $1 \leq j \leq m$. 

For each $j$, a subset of  vertex groups
in $\Theta_j$ are set stabilizers of domains in the orbit of the domains $V_j$, $W^{j,f}_t$. 
For every finite set of elements, $Fi$, in a vertex group that is not stabilized by
one of the set stabilizers, $W^{j,f}_t$, there is some constant $b_{Fi}$, 
such that for every quasimorphism from the subsequence, there exists a point $p_{Fi}$  in  some domain in the orbit of  $V_j$ (the point depends
on the quasimorphism), that is moved a distance bounded by no more than $b_{Fi}$ by the images (under the quasimorphisms $\{\psi^j_s\}$) 
of all the elements in the finite set $F_i$.

Furthermore, each set stabilizer, $W^{j,f}_t$, is generated by the (finitely many) edge groups that are connected to the vertex that it stabilizes, $\{E^{j,e}\}$, 
together with finitely many elements.
\endproclaim

\nfp By our constructions, the basepoints $\{c^j_s\}$, are displaced by a uniformly bounded distance in the projection complex, $P_K^j$, by the elements
$\psi^j_s(h_i)$, $1 \leq i \leq \ell$. 
Hence, 
for each quasimorphism, $\psi^j_s$, the segment: 
$[c_s^j,\psi^j_s(h_i)(c^j_s)]$, $1 \leq i \leq \ell$,
 is supported on boundedly many (embeddings of) domains in the orbit of $V_j$ in $C_K^j$.

Since the action of $H$ on the HHS $X$ is assumed to be cocompact, 
by conjugating the quasimorphisms $\{\psi^j_s\}$, we may assume that all the chosen basepoints, $\{c^j_s\}$, are some fixed point $p_0 \in C_K^j$.

By passing to a further subsequence of quasimorphisms, we gradually construct a simplicial bipartite tree $T$ on which the limit group $L^j_t$ acts
by isometries.
We start with the union of the collection of paths $[p_0,\psi^j_s(h_i)(p_0)]$, $1 \leq i \leq \ell$, in $C_K^j$. 

We denote on each of these paths points in
which the path moves from a subsegment of length bigger than $2K$ in a domain in the orbit of $V_j$ to a subsegment of length bigger than
$2K$ in another domain in the orbit of $V_j$ in $C_K^j$. Note that the number of points that we added along each such path is uniformly bounded.
By passing to a subsequence of quasimorphisms, we may assume that the combinatorial places of the points on
the union of these paths is identical for the entire subsequence.
   
If there exist points that we marked on the union of  paths $[p_0,\psi^j_s(h_i)(p_0)]$ for which there exists a subsequence 
such that these points remain in a bounded distance in 
$C_K^j$ under the images of the quasimorphisms in the entire subsequence, we pass to that subsequence, and consider them to be in the same equivalence class.

With the union of paths $[p_0,\psi^j_s(h_i)(p_0)]$ we further associate a finite combinatorial tree. The vertices are the equivalence classes of points
that we marked on the paths, and vertices for each of the domains that contain at least one stably unbounded subsegment. We connect a vertex that is 
associated with
a domain to the vertex that is associated with an equivalence class that contains a marked point in this domain.    

We continue by adding all the paths from $p_0$ to the images of elements (under the sequence of quasimorphisms) 
of length 2 in $H$, i.e., to images of elements of the form $h_{i_1}h_{i_2}$. We mark points, pass to subsequences, and define equivalence 
classes of points and vertices in exactly the same way as we did with the paths in the ball of radius 1.

We continue iteratively, each time looking at all the paths from $p_0$ to images under the sequence of quasimorphisms of elements of the ball of radius 
$n$ in the Cayley graph of $L^j_t$, w.r.t.\  the fixed set of generators. At each step we pass to a subsequence of the quasimorphisms, such that the
combinatorial position of the points that we mark, and the collection of stably bounded and stably unbounded segments will agree along the subsequence.

Finally, by taking a diagonal subsequence, we associate a simplicial tree, $T_j$, with the action of $L^j_t$ on $C_K^j$. The tree is bipartite, and contains vertices
that are associated with domains that contain  unbounded segments, and vertices that are associated with equivalence classes of points that
mark the transition between such domains. Since $L^j_t$ is f.g.\ and the segments $[p_0,\psi^j_s(h_i)(p_0)]$ are supported on boundedly many 
domains in $C^k_j$, the graph of groups decomposition that is associated with the action of $L^j_t$ of $T_j$ is bipartite and finite, and we denote it
$\Theta_j$.

By construction, all the vertex groups in $\Theta_j$ that are associated with equivalence classes of points that mark the 
transition from one domain to another
are stably bounded. Since the graph is bipartite,  all the edge groups in $\Theta_j$ are connected to such vertex groups. 
Hence, all the edge groups in $\Theta_j$, that we denote
$E^{j,e}$, are stably bounded subgroups.  


Since $L^j_t$ is f.g.\ each vertex group in $\Theta_j$ is generated by (finitely many) edge groups that are connected to it, in addition
to finitely many elements. In particular, this is true for the vertices that stabilize domains in $\Theta_j$, that we denote $W^{j,f}_t$, and the
conclusions of the proposition follow.

\line{\hss$\qed$}

In the sequel we study resolutions of the domain (set) stabilizers $W^{j,f}_t$ and not of the f.g.\  limit groups $L^j_t$. These domain stabilizers need
not be f.g.\ but proposition 3.5 proves that they are f.g.\ relative to finitely many stably bounded subgroups. 
This relative f.g.\ and the assumption that the domain stabilizers
act weakly acylindrically on the domains that they stabilize, suffice in order to apply the techniques of [Gu] and [Gr-Hu] to analyze the (superstable) actions
of the domain stabilizers $W^{j,f}_t$ on real trees, when these real trees are obtained as limits of their actions on domains in the orbits 
$HV_j$, $1 \leq j \leq m$.

\vglue 1.5pc
\proclaim{Lemma 3.6} With the assumptions of theorem 3.4 and proposition 3.5, at least one of the following holds:
\roster
\item"{(1)}"  there exists a subsequence of the quasimorphisms $\{\psi^j_s\}$ for which the traces of each element in $W^{j,f}_t$
are  bounded. In this case there exists a bounded subsequence, that we still denote, $\{\psi^j_s\}$, for which there exist points 
$p_s$, such that for every element $h \in H$, which is  in the preimage of $W^{j,f}_t$, there is a bound on the distances:
$d_{V_j}(p_s,\psi^j_s(h)(p_s))$. Note that the bounds depend on the element $h$, but not on the index $s$.

\item"{(2)}" the sequence $\{\psi^j_s\}$ contains a subsequence
(still denoted $\{\psi^j_s$\}), that converges into a non-trivial action of $W^{j,f}_t$ on some real tree $Y_{j,f}$ in which the subgroups $E^{j,e}$ are
elliptic.
\endroster
\endproclaim

\nfp The lemma is a natural generalization of Paulin's theorem [Pa], to countable groups that are generated by finitely many bounded subgroups in addition to
finitely many elements.

Suppose that there exists a subsequence of the quasimorphisms $\{\psi^j_s\}$ for which the traces of all the elements are bounded. We still denote
the subsequence $\{\psi^j_s\}$. By proposition 3.5,
a set stabilizer $W^{j,f}_t$ is generated by finitely many subgroups $E^{j,e}$ in addition to finitely many elements. Furthermore, for each subgroup $E^{j,e}$ 
there are points $p_{e,s}$, such that the images of each element in $E^{j,e}$, under the subsequence $\{\psi^j_s\}$, move them a bounded distance, where the
bound depends only on the specific element in the preimage of $E^{j,e}$ and not on the index $s$.

If there are no subgroups $E^{j,e}$, or if they are all f.g., $W^{j,f}_t$ is f.g.\ and the conclusion of part (1) follows by the same argument that given
an action of a f.g.\ group on a tree, 
if every element in a f.g.\ group has a fixed point, the whole group has a fixed point.

Suppose that there is
at least one non-f.g.\ subgroup $E^{j,e}$, that we denote $E^{j,1}$. 
Let $E^{j,1}$ and $E^{j,2}$ be two elliptic subgroups. There exist pairs of points: $(p_{s,1},p_{s,2})$ such that images of elements in the preimage of $E^{j,1}$ 
(under the sequence of quasimorphisms $\{\psi^j_s\}$) move the points 
$p_{s,1}$ a bounded amount (where the bounds depend on the element and not on the index $s$). Similarly,
the images of elements in the preimage of $E^{j,2}$ move the points $p_{s,2}$ a bounded amount.

If there exists a sequence of indices for which the distances between $p_{s,1}$ and $p_{s,2}$ are bounded, then we pass to this subsequence, and images
(under quasimorphisms in the subsequence) of elements in
the preimage of the subgroup that is generated by $E^{j,1}$ and $E^{j,2}$  move the points $p_{s,1}$ a bounded amount.  

Suppose that the distances between $p_{s,1}$ and $p_{s,2}$ are not bounded. We look at two sequences of f.g.\  subgroups,  $H^1_n$ and $H^2_n$, in the preimages
of $E^{j,1}$ and $E^{j,2}$ in correspondence, such that $H^i_n<H^i_{n+1}$, $i=1,2$, and the sequence of f.g.\ subgroups $H^1_n$ and $H^2_n$ approximate the
preimages (i.e., their corresponding unions are the entire preimages). We fix finite generating sets of each of the subgroups $H^1_n$ and $H^2_n$, and assume that
the fixed generating set of $H^i_n$ is a subset of the fixed generating set of $H^i_{n+1}$, for $i=1,2$ and every index $n$. 

For each index $s$ we look at a geodesic path between $p_{s,1}$ and $p_{s,2}$. For each pair of indices $n,s$, we define a subset $B^1_{n,s}$ to be a subset of the
geodesic path between $p_{s,1}$ and $p_{s,2}$, that contain all the points that move by the images of each element $h_1$ from the fixed set of generators of 
$H^1_n$, under the sequence of quasimorphisms $\{\psi^j_s\}$, a distance that is bounded by the bound on $d_{V_j}(p_{s,1},\psi^j_s(h_1)(p_{s,1})$ plus $10 \delta_j$
(where $\delta_j$ is the hyperbolicity constant of $V_j$). Similarly, we define $B^2_{n,s}$ as a subset of the geodesic between $p_{s,1}$ and $p_{s,2}$,
as sets of points with  similar bounds on the distances that they are moved by the fixed set of generators of the subgroups $H^2_n$.  

Note that $p_{s,i} \in B^i_{n,s}$, $i=1,2$.  
Since the fixed set of generators of $H^i_n$ is contained in the fixed set of generators of $H^i_{n+1}$, $B^i_{n+1,s} \subset B^i_{n,s}$.
If there exists some index $n_0$, for which the sequence of   distances between
$B^1_{n,s}$ and $B^2_{n,s}$ is not bounded, then there exists an element $h$ in the preimage of the subgroup that is generated by 
$E^{j,1}$ and $E^{j,2}$,  for which the traces of the images $\{\psi^j_s(h)\}$ are unbounded.  

Therefore, under the assumptions of part (1), for every index $n$, the distances between the sets $B^1_{n,s}$ and $B^2_{n,s}$ are bounded.
We fix an index $s$. If for some index $n_0$, $B^1_{n_0,s}$ and $B^2_{n_0,s}$ have empty intersection, we pick a point $p_{s,3}$ to be a point on the geodesic
from $p_{s,1}$ to $p_{s,2}$ that is not in the union of $B^i_{n_0,s}$. 

If for every index $n$, $B^1_{n,s}$ intersects non-trivially $B^2_{n,s}$, we pick $p_{s,3}$ to be a point in the
(non-empty) intersection of $\cap B^1_{n,s}$ and $\cap B^2_{n,s}$.

Since the distances between $B^1_{n,s}$ and $B^2_{n,s}$ are bounded, and
$B^i_{n+1,s} \subset B^i_{n,s}$, the images of every element $h$ in the subgroup that is generated by $E^{j,1}$ and $E^{j,2}$ under the sequence
of quasimorphisms $\{\psi^j_s\}$, move $p_{s,3}$ a bounded distance, where the bound depends on the element $h$ and not on the index $s$.

\smallskip
So far we proved that under the assumptions of part (1) of the lemma, for every pair of elliptic subgroups, $E^{j,e_1}$ and $E^{j,e_2}$, there exists a point
$p_s$ that moves a bounded distance by the images of every given element in the preimages of the subgroup that is generated by the two elliptic subgroups.
Continuing inductively, there exists a point (still denoted $p_s$) that is moved a bounded distance by the images of every given element in the preimage
of the subgroup that is generated by all the (finitely many) elliptic subgroups $E^{j,e}$.

By proposition 3.5, the set stabilizer $W^{j,f}_t$ is generated by finitely many stably bounded subgroups, $E^{j,e}$, in addition to finitely many elements.
By the assumptions of part (1) there is a bound on the traces of each fixed element under the sequence of quasimorphisms $\{\psi^j_s\}$. Hence, the subgroup
that is generated by any given element in $W^{j,f}_t$ is stably bounded as well. Therefore, $W^{j,f}_t$ is generated by finitely many
stably bounded subgroups, so for each $s$ and every element $h$ in the preimage of $W^{j,f}_t$, there is a bound on the distance that the images 
$\{\psi^j_s(h)\}$ move some chosen basepoint $p_s$, and the conclusion of part (1) follows.

To prove part (2) suppose that there exists an element $u$ in the preimage of $W^{j,f}_t$, for which the sequence of traces of the elements 
$\{\psi^j_s(u)\}$ is unbounded. In this case we pass to a subsequence for which the sequence of traces of the elements $\{\psi^j_s(u)\}$ does not have a 
bounded subsequence, and add the image of the element $u$ to the generating set of $W^{j,f}_t$. 

The set stabilizer $W^{j,f}_t$ is obtained as a limit from a sequence of quasimorphisms. We set $N^{j,f}_t$ to be the normal subgroup of 
$W^{j,f}_t$ that act stably quasi-trivially on the domain $V_j$, and set $Q^{j,f}_t$ to be the quotient:
$W^{j,f}_t/N^{j,f}_t$. In studying a limit action of $W^{j,f}_t$, we are actually studying the limit action of the (limit) quotient $Q^{j,f}_t$. 

$Q^{j,f}_t$, the quotient group,  is generated by finitely many 
infinitely generated bounded subgroups, the images (under the quotient map) of some bounded subgroups $E^{j,e}$, in addition
to finitely many elements, and the image
of the element $u$, for which the sequence of traces of $\{\psi^j_s(u)\}$ has no bounded subsequence.

Note that if there are no bounded (edge)  subgroups, or if the images of all the bounded subgroups, $E^{j,e}$, in $Q^{j,f}_t$ are f.g., 
then the conclusion of part (2)
of the lemma follows from Paulin's theorem [Pa]. Hence, we assume that there is at least one bounded (edge) subgroup, $E^{j,e}$,  
with a non-f.g.\ image in
$Q^{j,f}_t$.
 
For each index $s$, let $A_{s}(u)$ be an axis of the loxodromic element $\psi^j_s(u)$. Let $p_{s,e}$ be points that are moved a bounded distance by
the image of each element in the preimage of $E^{j,e}$, under the sequence of quasimorphisms $\{\psi^j_s\}$. Given a point $p_{s_e}$, we set $b_{s,e}$ to
be one of the closest points to $p_{s,e}$ in $A_{s}(u)$. 

For each index $s$, we look at the convex hull of the points $b_{s,e}$ in the axis $A_s(u)$. We set $b_s$ to be a point in a bounded distance from the middle
of the convex hull in $A_s(u)$ of the points, $b_{s,e}$. $b_s$ will serve as a basepoint for the iterative construction of the limit tree.

The bounded subgroups $E^{j,e}$, with non-f.g.\ image in $Q^{j,f}_t$, are countable and so are their preimages. We fix a sequence of
f.g.\ approximations for the preimages of the bounded subgroups $E^{j,e}$, that we denote $\{H^e_n\}$. We assume that for each $e$, $H^e_n < H^e_{n+1}$.
Furthermore, we fix finite generating sets of each of the approximating subgroups $H^e_n$, and assume that the finite generating set of $H^e_n$ is contained
in the finite generating set of $H^e_{n+1}$.

We start with the finite generating sets of the subgroups $H^e_1$, and the additional finitely many generators, that include the element $u$. We look
at the sequence of images of these (finitely many) elements under the sequence of quasimorphisms, $\{\psi^j_s\}$, and the image of the basepoints
$\{b_s\}$ under these images. By rescaling so that the maximal distance between $b_s$ and its images will have length $1$, we can pass to a subsequence that
converges into a finite tree , of diameter at most 2. 

We gradually enlarge the finite generating sets, to include the (prefixed) finite generating sets 
of the subgroups, $H^e_n$. For each $n$ we further rescale the metric so that the base points $\{b_s\}$ move a maximal distance $1$ by the images of the
enlarged finite set of generators, and pass to a subsequence that converges into a finite tree, of diameter at most 2. Note that the tree
that was obtained from the sequence at step $n$ is embedded by homothety into the tree that is obtained in step $n+1$. Also, note that because of our 
weak acylindricity
assumption the homothety constant of embeddings between consecutive trees can be strictly smaller than 1 only for boundedly many indices $n$.   

Finally, we take a diagonal sequence of quasimorphisms, and obtain a tree, $Y_1$, of diameter at most 2, in which the (infinite) sets of generators of
the images of each of the elliptic subgroups, $E^{j,e}$, in $Q^{j,f}_t$ fix points.

Now we gradually add elements that can be presented as words of length 2. First as words of length 2 in the (prefixed) generating set of $H^e_1$ and 
the additional finitely many generators, and then gradually adding elements that can be presented as words of length 2 in the prefixed generating
sets of $H^e_n$ and the additional finitely many generating sets. Note that the rescaling of the metric was already done in constructing the tree $Y_1$.
By passing to another diagonal sequence of quasimorphisms, we get convergence into
a tree $Y_2$, of diameter at most 4.

We continue iteratively by adding elements that can be presented as words of larger and larger length, and eventually obtain a (diagonal) subsequence
of quasimorphisms that converges into an action of the group $Q^{j,f}_t$ on some real tree $Y_{j,f}$. By construction, the images of the bounded
subgroups $E^{j,e}$ in $Q^{j,f}_t$ fix points in $Y_{j,f}$. Finally, the action of $Q^{j,f}_t$ on $Y_{j,f}$ is non-trivial, since either the image of the base
point in $Y_{j,f}$ is on the axis of the image of $u$ which is loxodromic in $Q^{j,f}_t$, or the base point is fixed by the infinite cyclic subgroup 
that is generated by the
image of $u$, and the base point is not fixed by either one of the elliptic subgroups, $E^{j,e}$, or by one of the additional (finitely many) generators. 

\line{\hss$\qed$}

The actions of the quotients of the set stabilizers $Q^{j,f}_t$ on the trees $Y_{j,f}$, that are constructed in lemma 3.6,
 enable one to apply the techniques of [Gu] and [Gr-Hu] to analyze these actions.

\vglue 1.5pc
\proclaim{Lemma 3.7 (cf. [Gu], Main Theorem)} Suppose that part (2) in lemma 3.6 holds, and $Q^{j,f}_t$ acts on some limit tree $Y_{j,f}$.
Then either $Q^{j,f}_t$ splits over the stabilizer of a tripod in $Y_{j,f}$, which is a finite (uniformly bounded) subgroup, 
or over the stabilizer of an
unstable segment in $Y_{j,f}$, that is a  finite (bounded) subgroups as well,
splittings in which all the subgroups
$E^{j,e}$ are elliptic, 
or $Y_{j,f}$ has a decomposition into a graph 
of actions where each vertex action is either simplicial, axial, or of IET type (see [Gu] for these notions).
\endproclaim

\nfp  By lemma 3.5 the groups $Q^{j,f}$ are generated by finitely many bounded subgroups that fix points in the trees
$Y_{j,f}$, in addition to finitely many elements. Since we assumed that the set stabilizer of a domain $V_j$ acts weakly
acylindrically on $V_j$, the stabilizer of an unstable segment and the stabilizer of a tripod  in $Y_{j,f}$ is finite and universally bounded (see lemma 4.7
in [Gr-Hu]). 

Hence, the actions of $Q^{j,f}$ on $Y_{j,f}$ satisfy the assumptions of the main theorem (Theorem 5.1) in [Gu]. The conclusion of this 
theorem is the conclusion of the lemma.  

\line{\hss$\qed$}

By proposition 3.5, for each $j$, $1 \leq j \leq m$, the set stabilizers, $W^{j,f}_t$, are generated by finitely many bounded subgroups $E^{j,e}$ that they include, 
together with finitely many elements. This enables one to adopt the techniques
of Weidmann [We1], and obtain an acylindrical accessibility principle with a family of graph of groups decompositions of the quotients of the set stabilizers,
 $Q^{j,f}_t$.

\vglue 1.5pc
\proclaim{Lemma 3.8} With the notation of lemma 3.5, 
recall that $Q^{j,f}_t$ is generated by the images of finitely many
bounded subgroups, $E^{j,e}$, together with finitely many elements.
$Q^{j,f}_t$, the quotient of the domain stabilizer $W^{j,f}_t$ by the subgroup $N^{j,f}_t$  that act quasi-trivially on the domain,
 satisfies a relative acylindrical accessibility principle.

For every two integers, $(k,c)$,  there exists some bound $b_{k,c}$, such that in every $(k,c)$ acylindrical graph of groups decomposition 
of $Q^{j,f}_t$ in which
the images of the bounded subgroups, 
$\{E^{j,e}_t\}$, that are contained in $W^{j,f}_t$, are elliptic, there are  no more than $b_{k,c}$ edges. 
\endproclaim

\nfp This is a relative version of [We1], and follows by the same arguments according to [We2]. It appears as proposition B.3 in [GHL].

\line{\hss$\qed$}

Lemmas 3.7 and 3.8  enable us to associate a  graph of groups decomposition with the action of $W^{j,f}_t$ on $V_j$, or rather 
a graph of groups decomposition of the quotient $Q^{j,f}_t$ that is associated with its action on the  limit tree $Y_{j,f}$.

\vglue 1.5pc
\proclaim{Proposition 3.9} With the action of $Q^{j,f}_t$  on the real tree $Y_{j,f}$ it is possible to canonically 
associate a graph of groups 
decomposition $\Lambda_{j,f}$. The vertex groups in this graph of groups are point stabilizers, and stabilizers of axial and IET components. The edge groups
are either uniformly finite groups, or stabilizers of stable segments   and cyclic extensions of stabilizers of tripods in $Y_{j,f}$. 
\endproclaim

\nfp We start with the action of $Q^{j,f}$ on the real tree $Y_{j,f}$ according to lemma 3.7. 
If $Q^{j,f}$ splits over a stabilizer of a tripod, or over the
stabilizer of an unstable segment, then  $Q^{j,f}$ splits over a finite group. We denote this splitting $\Theta_1$.

Since we assumed that domain stabilizers act weakly acylindrically on the domain that they stabilize, the orders of the finite groups 
over which $Q^{j,f}_t$ splits, i.e., the edge groups
in $\Theta_1$, are uniformly bounded (cf. lemma 4.7 in [Gr-Hu]).
Furthermore, since the images of the bounded subgroups $E^{j,e}$ fix points in $Y_{j,f}$, 
the images of the subgroups $E^{j,e}$ are all elliptic in $\Theta_1$. 

Since $Q^{j,f}_t$ is generated by the images of the (finitely many) elliptic subgroups
$E^{j,e}$ in addition to finitely many elements, and the edge groups in $\Theta_1$ are finite, each vertex group in $\Theta_1$ is generated by finitely
many bounded subgroups in addition to finitely many elements. 

Now we can restrict the sequence of quasimorphisms to each vertex group in $\Theta_1$, pass to a subsequence, and apply  lemmas 4.6 and 4.7 to the vertex
subgroups. Suppose that at least  one of the vertex groups splits over a finite group, that has to be of uniformly bounded order that depends only on
the acylindricity constants of the action of the domain stabilizers in $H$. In that case we can further refine $\Theta_1$,
and obtain another splitting  of $Q^{j,f}$ with finite (uniformly bounded) edge groups with strictly more edges. By the acylindrical accessibility principle
that is stated in lemma 3.8, this refinement procedure terminates after finitely many steps.

When the successive refinement of splittings along bounded edge groups terminate, by lemma 3.7 every vertex group in the constructed splitting
obtains a (possibly trivial) splitting in which vertex groups are point stabilizers, virtually abelian and QH vertex groups. Edge groups 
are all virtually abelian. By construction, the images of the bounded subgroups $E^{j,e}$ must be elliptic in the constructed splitting. 

\line{\hss$\qed$}

Finally, starting with a subsequence of quasimorphisms $\{\psi^j_s\}$ and applying propositions 3.5 and 3.9
to the associated quasi-actions of $H$ on the quasi-tree $C_K^j$,
 it is possible to modify the techniques that were used in the proof of theorem 2.7, and pass to a convergent subsequence from which it is possible 
to construct a resolution for each of the quotient groups $Q^{j,f}_t$ (which are quotients of the domain stabilizers $W^{j,f}_t$), 
that terminate with a limit group with  bounds on the basepoint displacements (in $C_K^j$) 
of its associated quasimorphisms.

\vglue 1.5pc
\proclaim{Proposition 3.10} Let $\{\psi^j_s\}$ be a sequence of quasimorphisms, that were obtained from a sequence of automorphisms $\{\varphi_s\}$ 
using a resolution that was constructed from twisted actions of $H$ on the projection complex $P_K^j$. Recall  that $\{\psi_s^j\}$ converges into the
 limit group, $L^j_t$, and there  are global bounds on the displacements of the images of the fixed finite set of generators of $L^j_t$ under the 
quasimorphisms when acting on the projection complex $P_K^j$. 

Then it is possible to pass to a subsequence of the quasimorphisms, and when viewed as quasimorphisms into the isometry group of the
quasi-tree of metric spaces, $C_K^j$, it is possible to associate with the subsequence finitely many  (finite) resolutions of the 
quotients, $Q^{j,f}_t$, of the set stabilizers
of domains in the orbit of the domain $V_j$, $W^{j,f}_t$.

Furthermore, the resolutions terminate with a terminal limit group, and with a sequence of quasimorphisms that converges into it, $\{\eta^{j,f}_s\}$. 
With the  terminal limit group there is an associated finite graph of groups decomposition, in which all the edge groups are finite,
and all the vertex groups are stably bounded subgroups.  
\endproclaim

\nfp  The proof is essentially identical to the argument that was used in the proof of theorem 2.7 and in section 6 of [Gr-Hu] in the f.g.\ case.

A set stabilizer of a domain in the orbit of the domain $V_j$, $W^{j,f}_t$, is a subgroup of the f.g.\ limit group $L^j_t$. 
We fix a finite generating set of $L^j_t$. This allows us
to partially order the elements of $L^j_t$ according to their length in the Cayley graph of $L^j_t$, where each equivalence class contains the elements
of fixed length in $L^j_t$. Hence, each equivalence class is finite.
$W^{j,f}_t<L^j_t$, so $W^{j,f}_t$ inherits this partial order. Each equivalence class in $W^{j,f}_t$ is finite or empty. 

We say that a quotient $LW$ of $W^j_t$ is a limit group if it is obtained as a limit of a sequence of quasimorphisms of some infinitely generated free group 
$F_{\infty}$ into the 
isometry group of the domain $V_j$, such that the map: $F_{\infty} \to LW$ factors as: $F_{\infty} \to W^{j,f}_t \to Q^{j,f}_t \to LW$. 

Let $\nu^j_s$ be a convergent sequence of quasimorphisms of $F_{\infty}$ into the isometry group
of $V_j$, and let its limit (group) be $LW$. As in the proof of theorem 2.7, we say that an element $u \in F_{\infty}$ is stably bounded, 
if the traces of the images 
of $u$ under the sequence of quasimorphisms are bounded. It is stably unbounded if the sequence of traces has no bounded subsequence. The image of $u$ in $LW$ is 
called stably bounded or stably unbounded, if some preimage of it in $F_{\infty}$ is.

As in the proof of theorem 2.7, by passing to a subsequence, we can assume that every 
element in $LW$ (and hence in $F_{\infty}$ and in $W^{j,f}_t$ and $Q^{j,f}_t$), is either stably bounded or stably unbounded 
(w.r.t. the subsequence of quasimorphisms). 
In the rest of
the proof we will consider
only convergent sequences of quasimorphisms (and their corresponding limit groups) for which the images of elements in the subgroups $E^{j,e}_t<W^{j,f}_t$, 
that were defined in proposition 3.5, are stably bounded.

We continue by adapting the argument that appears in section 1 in [Ja-Se] and in section 6 in [Gr-Hu]. Suppose that there exists an infinite 
(properly) descending chain of limit groups:
$$Q^{j,f}_t=LW_0 \, \to \, LW_1 \, \to \, LW_2 \, \to \, LW_3 \, \to \ldots$$
of limit quotients of $Q^{j,f}_t$, in which the images of the subgroups $E^{j,e}_t$ are all stably bounded, and non-trivial stably bounded elements in $LW_d$ 
are mapped
to non-trivial stably bounded elements in $LW_{d+1}$. We further assume that if an element $u \in F_{\infty}$ is stably 
bounded w.r.t.\ to the sequence of quasimorphisms
that converges to $LW_d$, then the global bound on the traces of the images of $u$ under the sequence of quasimorphisms that converge into
$LW_d$, is also a global bound on the traces of the images of $u$ by the sequence of quasimorphisms that converge into $LW_{d+1}$.

We fix $LW_1$ as a proper limit quotient of $Q^{j,f}_t$, in which the non-trivial stably bounded elements in $Q^{j,f}_t$ are mapped to non-trivial 
stably bounded
elements, and $LW_1$ can be extended to an infinite proper descending chain of limit quotients of itself. We further require
that global bounds on traces of images of stably bounded 
elements in $F_{\infty}$ under the quasimorphisms that converge into $Q^{j,f}_t$ remain global bounds for traces of the same elements 
under the sequence of quasimorphisms that converge into 
$LW_1$. We pick $LW_1$ to have the maximum possible elements of length
1 (if there are such in $W^j_t$) that are mapped to the trivial element among such limit quotients of $Q^{j,f}_t$.

We continue iteratively as in [Ja-Se] and [Gr-Hu]. We set $LW_{d+1}$ to be a proper limit quotient of $LW_d$,  such that the non-trivial 
stably bounded elements in 
$LW_{d}$ are mapped to non-trivial stably bounded elements in $LW_{d+1}$, and $LW_{d+1}$ has an infinite descending chain that satisfies the iterative properties.
We further require
that global bounds on traces of images of stably bounded 
elements in $F_{\infty}$ under the quasimorphisms that converge into $LW_d$ remain global bounds for traces of the same elements 
under the sequence of quasimorphisms that converge into 
$LW_{d+1}$. We pick $LW_{d+1}$ to have the maximum possible elements of length $d+1$
that are mapped to the identity, among all such limit quotients of
$LW_d$.

We obtained an infinite sequence as above. $F_{\infty}$ is a countable group, and we order its elements. For each index $d$, we pick a quasimorphism of
$F_{\infty}$ to $Isom(C_K^j)$, $\tau_d$,  from the sequence that converges to $LW_d$, 
that maps the first $d$ elements to stably trivial, non-trivial, and stably bounded and unbounded elements, according to
their image in $LW_d$. If an element among the first $d$ elements is stably unbounded, we require that the quasimorphism $\tau_d$ will send it to an element
with trace  bigger than $d$.

The sequence $\{\tau_d\}$ subconverges to a limit group $LW_{\infty}$, in which the images of the subgroups $E^{j,e}$ are stably bounded, 
and $LW_{\infty}$ is generated by these 
(bounded) subgroups in addition to finitely many elements (since it is a quotient of $W^j_t$). By construction $LW_{\infty}$ is the direct limit of
the sequence of limit groups: $LW_1 \to LW_2 \to \ldots$.

Following section 1 of [Ja-Se] and section 6 of [Gr-Hu], by further passing to a subsequence 
of the quasimorphisms, $\{\tau_d\}$, we obtain a finite resolution:
$$LW_{\infty} \, \to \, LW^1_{\infty} \, \to \, \ldots \, \to \, LW^r_{\infty}$$
where $LW^r_{\infty}$ is associated with  a graph of groups decomposition in which all the edge groups are finite,
and the elements in all the  vertex groups are stably bounded.

With each of the limit groups $LW^c_{\infty}$ there is an associated (finite) graph of groups decomposition with virtually abelian edge groups in which
the images of the subgroups $E^{j,e}$ are elliptic (since they are all stably bounded).
Since the actions of the set stabilizers in $H$ on $V_j$ are weakly acylindrical, and the elements in the vertex groups in the virtually abelian decomposition
that is associated with the terminal limit group $LW^r_{\infty}$ are all stably bounded, the edge groups in all the graphs of groups decompositions
that are associated with the limit groups $LW^c_{\infty}$ are f.g.\ virtually abelian.

With $LW_{\infty}$ there is an associated finite graph of groups decomposition with f.g.\ virtually abelian edge groups. Each vertex group in this 
decomposition is generated by finitely many bounded subgroups together with finitely many elements. Hence, each vertex group
 inherits a finite decomposition that is inherited from the finite graph of groups decomposition that is associated with $LW^1_{\infty}$.
We continue iteratively with vertex groups and the virtually abelian decompositions that are associated with the limit groups according to the
finite resolution. Note that the edge groups in all these graphs of groups decompositions are f.g.\ virtually abelian.

The images of the subgroups $E^{j,e}$ being bounded and embedded in $LW_{\infty}$, remain elliptic in this finite iterative sequence of virtually
abelian graphs of groups decompositions. Since $LW_{\infty}$ is generated by the images of the bounded subgroups $E^{j,e}$ in addition to finitely
many elements, and all the edge groups in the virtually abelian decompositions of $LW^c_{\infty}$ are f.p.\ with all subgroups f.p.,
$LW_{\infty}$ is generated by finitely many stably bounded subgroups in addition to finitely many elements and finitely many relations. Furthermore,
each of these finitely many stably bounded subgroups, is generated by finitely many conjugates of images of some of the bounded subgroups $E^{j,e}$, together with
finitely many elements.

Hence, $LW_{\infty}$ is generated by finitely many stably bounded subgroups, that are images (in fact isomorphic images) of bounded subgroups is some 
limit group in the first sequence, $LW_{d_1}$, in addition to a finite collection of elements and a finite collection of relations.
Hence,
it follows that for some larger $d_2$ and all $d>d_2$, $LW_d$ is a quotient of $LW_{\infty}$.
This contradicts the assumptions that the sequence $LW_1 \to LW_2 \to \ldots$ is a sequence of proper quotients  in which stably bounded subgroups are embedded,
and $LW_{\infty}$ is the direct limit of the sequence. 

Therefore, every sequence of proper quotients of limit quotients of $Q^{j,f}_t$, $LW_1 \to LW_2 \to \ldots$,
in which:
\roster
\item"{(i)}"
 the images of the subgroups $E^{j,e}$ are bounded. 

\item"{(ii)}"  a stably bounded 
subgroup in $LW_d$  is mapped isomorphically into a stably bounded subgroup in $LW_{d+1}$.

\item"{(iii)}"  a bound on the displacements of  the images of a stably bounded element under the sequence of quasimorphisms that converge into  $LW_d$, remains
a bound on the displacement of the images of that element under the sequence of quasimorphisms that converge into
$LW_{d+1}$.
\endroster
terminates after finitely many steps.  

Now, let $\{\psi^j_t\}$ be a sequence of quasimorphisms of $H$ into the isometry group of $P_K^j$ that converges into the limit group $L^j_t$, such that
there exists
a global bound on the displacements of the images of any given element in a vertex group in the graph of groups that is associated with the terminal
limit group $L^j_t$, under the sequence $\{\psi^j_t\}$. We look at $\{\psi^j_t\}$ as a sequence of
quasimorphisms into the quasitree of metric spaces $C_K^j$. 

By passing to the vertex groups in the terminal graph of groups decomposition of $L^j_t$,
and to a convergent subsequence, we define the bounded subgroups $E^{j,e}$ w.r.t. the
subsequence of quasimorphisms according to proposition 3.5. We further pass to a subsequence of the sequence of quasimorphisms,
and assume that with each set stabilizer of a domain in the orbit of 
$V_j$ (that are defined in proposition 3.5), $W^{j,f}_t$, it is possible to associate a limit action on some real tree $Y_{j,f}$ (see lemma 3.7).
By construction the subgroups $E^{j,e}$ in $W^{j,f}_t$, that were defined in proposition 3.5, are bounded and, hence, fix points in these limit trees $Y_{j,f}$.
Furthermore, by proposition 3.9, $W^{j,f}_t$ inherits a finite graph of groups decomposition from its action on $Y_{j,f}$.

Let $Q^{j,f}_t$ be the quotient of $W^{j,f}_t$ by the stable  quasi point stabilizer of $V_j$. $Q^{j,f}_t$ 
inherits a virtually abelian decomposition 
from the finite graph of groups decomposition of $W^{j,f}_t$. With the virtually abelian decomposition of $Q^{j,f}_t$ we can associate modular groups according
to section 5.4 in [Gr-Hu]. We use these modular groups to shorten the quasimorphisms $\{\psi^j_t\}$, and pass to a convergent sequence of quasimorphisms,
that converges into a quotient of $Q^{j,f}_t$. 

If the quotient is isomorphic to $Q^{j,f}_t$ it is possible to further refine the virtually abelian decomposition of $Q^{j,f}_t$. 
By lemma 3.7 such a refinement can occur
only finitely many times. Hence, after finitely many such shortenings we get  either a proper quotient of $Q^{j,f}_t$, or we get to the terminal step
of the resolution, in which all the edge groups are finite and all the vertex groups are stably bounded. 

By our previous argument we can pass from a limit group to a proper quotient of it only finitely many times. Hence, after finitely many steps
we must get to the final step of the resolution in which all the vertex groups are stably bounded, and all the edge groups are (uniformly bounded) finite. 

\line{\hss$\qed$}

By propositions 3.5 and 3.10, given a sequence of automorphisms, $\{\varphi_s\}$ in $Aut(G)$, it is possible to pass to
a subsequence that converges into an $m$-collection of hybrid resolutions. Each hybrid resolution starts with a resolution of some limit quotient of
the finite index subgroup $H<G$, and continues with finitely many resolutions of the relative f.g.\  domain (set) stabilizers in the terminal limit group of the 
resolution we started with. 

As in theorem 2.7, that deals with the case of a product space, our next goal is to use the $m$-collections of hybrid resolutions 
that are obtained from convergent sequences of automorphisms, 
to construct a higher rank MR diagram, i.e.,
to obtain finitely many $m$-collections of hybrid cover resolutions, such that every automorphism in $Aut(G)$ factors through at least one of them. i.e.,
each automorphism factors through all the $m$ hybrid resolutions in one of the finitely many $m$-collections, where factors is in the sense of definition 2.6.

To find such finitely many $m$-collections of hybrid cover resolutions, we use a compactness argument, similar to the one that was used in
proving theorem 2.7 (in the product case). In hybrid resolutions, some of the groups
along the resolutions are f.g.\ and not f.p.\ and even worse, the limit subgroups that stabilize setwise a given domain are in general
not f.g.\ nor are the elliptic vertex groups that are connected to these domain stabilizers in the graphs of groups $\Theta_j$, 
 nor the edge groups that connect the stabilizers of the domains
to the elliptic vertex groups. Furthermore, the finitely many resolutions of the stabilizers of domains in each hybrid resolution,
are resolutions of subgroups that are in general not f.g.\ - they are f.g.\ relative to finitely many elliptic subgroups. 
Hence, we need to find $covers$ of hybrid resolutions that are f.p.\ objects, and that encode the entire geometry of the hybrid resolutions. 

Let $HbRes$
be a hybrid resolution. By construction, the top resolution in a hybrid resolution that we denote $TPRes$, i.e., the resolution that was constructed 
from actions of $H$ on
one of the  projection complexes, $P_K^j$, $1 \leq j \leq m$, is a resolution that has the same structure as the resolutions that were constructed in 
the case of product
spaces. Hence, with the resolution $TPRes$ we associate a cover in the same way that it was defined in [Ja-Se] and in definition 2.4.  i.e.,
the terminal limit group of the resolution $TPRes$ is replaced by some f.p.\ approximation, and on top of this terminal f.p.\ group we construct 
a completion into which the covers of all the limit groups along the cover resolution are embedded (see definition 2.4). Clearly, the whole 
completion is a f.p.\ group, that terminates in a f.p.\ limit group that we denote $L_t$.

The terminal limit group $L_t$ is equipped with a graph of groups decomposition along (uniformly bounded) finite groups, and we continue to the second
part of the hybrid resolution with each vertex group separately, in parallel. For brevity we denote each of the vertex groups in the graph of groups
with finite edge groups, $L_t$, which is also a f.p.\ group.
 
In proposition 3.5, we have associated a graph of groups decomposition $\Theta$, with each such vertex group (that is covered by $L_t$),
and $\Theta$ is part of the hybrid resolution. The graph of groups $\Theta$ is a bipartite graph of groups,
where some of the vertices are stably bounded, 
and the others are limits of set stabilizers of domains. Note that vertex and edge groups
in $\Theta$ are not necessarily finitely generated. Let $D_t$ be the fundamental group of $\Theta$, which is a f.g.\ (vertex) subgroup of 
the terminal limit group of the top resolution
$TPRes$ in the hybrid resolution $HbRes$. $L_t$ is a f.p.\ cover of $D_t$.
Hence, there is a finite set of elements from the various vertex groups in $\Theta$, that
together with the finite set of Bass-Serre generators in $\Theta$, generate $D_t$. We fix this generating set of $D_t$.

Let $W^f$ be a vertex group which is a limit of sets stabilizers of a domain in $\Theta$. By Proposition 4.5, $W^f$ is generated by finitely many stably bounded
subgroups 
(the edge groups that are connected to the vertex that is stabilized by $W^f$  in the graph of groups $\Theta$), and additional finitely many elements. 
By proposition 3.10 $W^f$ admits a finite resolution:  
$$W^f/K^f=W^f_1/K^f_1 \, \to \, \ldots \, \to \, W^f_r/K^f_r$$
where $K^f_i$, $i=1,\ldots,r$, are the stable kernels, i.e., the collection of elements that act stably quasi-trivially on the associated domain 
under the corresponding sequence of quasimorphisms.

By construction, the terminal limit group $W^f_r/K^f_r$ is equipped with a graph of groups decompositions in which all the edge groups
are finite and uniformly bounded, and all the vertex groups in this graph of groups are stably bounded. This implies that
for each finite set of elements $Fi$ in the preimage of such a vertex group in $H$, there is some bound $b_{Fi}$, such that for any quasimorphism
from the sequence that converges to $W^f_r$, there exists a point in the domain that is moved a distance bounded by $b_{Fi}$, by all the elements in
the set $Fi$.

$W^f$ is not necessarily a f.g.\ group, but since it is generated by finitely many stably bounded subgroups in addition to finitely many elements,  
there exists a f.g.\ subgroup $R<W^f$, that we can assume contains the finitely many elements in $W^f$ that are part of the fixed generating set of $D_t$, such that
$R$ obtains a resolution: $$R/KR=R_1/KR_1 \, \to \, \ldots \, \to \, R_r/KR_r$$
where $R_i<W^f_i$, $KR_i=K_i \cap R_i$, the virtually abelian decompositions that are associated with the groups $R_i/KR_i$ have the same structure as those of $W^f_i/K_i$,
and the quotient maps in the resolutions of $R/KR$ and $W^f/K$ together with the embedding of the groups $R_i$ into $W^f_i$ form a commutative diagram.

Altogether, we have a finite set of generators of $D_t$, the fundamental group of the graph of groups $\Theta$, 
to which we add finite sets of generators of the groups $R^f$ that approximate the limit stabilizers
of domains $W^f$ in $\Theta$ (and have resolutions with the same structure as those of $W^f$).
We continue with this fixed finite set of generators of $D_t$, that are all elements in vertex groups in $\Theta$ and Bass-Serre generators in $\Theta$.

Recall that  $D_t$ is a quotient of its f.p.\ approximation $L_t$. We further define an abstract  f.p.\ group, $U_t$, that is 
generated by copies of the fixed set generators of $D_t$, is naturally a quotient of $L_t$ and $D_t$ is a quotient of it, and it is
the fundamental group of a graph of  groups with f.g.\ vertex and edge groups that is similar to $\Theta$. i.e., the quotient map from $U_t$ onto $D_t$
preserve the graphs of group structures of both - vertex and edge groups are mapped into vertex and edge groups.
 
Note that the relations between the fixed finite set of elements in $D_t$, that generate $D_t$ and are all contained in vertex groups in $\Theta$ or are
Bass-Serre generators, can be taken to be relations in
each of the vertex groups, and relations that correspond to foldings. i.e., relations that correspond to an enlargement of an edge group. Since these relations
preserve the graph of groups structure, we can construct $U_t$ by adding finitely many such relations that will enable us to define a natural map from $L_t$ onto $U_t$,
and will guarantee that $U_t$ is the fundamental group of a graph of groups that has a similar graph of groups structure as $\Theta$, 
just that its vertex and edge groups are all
finitely generated (and $U_t$ is finitely presented). 
Furthermore,  the map from $U_t$ onto $D_t$ respects the graph of groups structure, i.e., vertex and edge groups in the graph of groups
decomposition of $U_t$ are mapped to
vertex and edge groups in $\Theta$.

At this point we are ready to define a $cover$ of the hybrid resolution, $HbRes$. 

\proclaim{Definition 3.11} 
The cover $CHbRes$ of $HbRes$ is constructed from:

\roster
\item"{(1)}"  a f.p.\ cover  of the 
completion of the top resolution, $TPRes$, that terminates with the f.p.\ group $L_t$.

\item"{(2)}" the f.p.\ group $U_t$.

\item"{(3)}" f.p.\ covers of the completions of the resolutions of the subgroups $R^f$, that are themselves f.g.\ approximations of the resolutions of
the limit set stabilizers of the domains, $W^f$ (that are vertex groups in the graph of groups $\Theta$.
\endroster.

$CHbRes$, the cover of the hybrid resolution $HbRes$, is the group that is generated by the f.p.\ covers from parts (1) and (2), and the f.p.\
group $U_t$, to which we add finitely many relations. We identify  the fixed set of  generators of $L_t$ (the terminal group in the cover
of $TPRes$) with their image in $U_t$. We further identify the elements that are associated with the generators of each of the group $R^f$ in $U_t$,
with the elements that are associated with these generators in the f.p.\ cover of the corresponding resolution of the group $R^f$ from part (3).

Clearly, the group that is associated with $CHbRes$ is f.p.\ and the convergent subsequence of automorphisms $\{\varphi_s\}$ that was used to obtain the
hybrid resolution, $HbRes$, asymptotically factors through it. i.e., all the automorphisms from the convergent subsequence factor through $CHbRes$, except 
for at most finitely many of them.
\endproclaim

Once we defined a (f.p.) cover of a hybrid resolution, we can apply a compactness argument similar to the one that was used to prove
theorem 2.7, and obtain  finitely many $m$-collections of covers of hybrid resolutions,
that form a higher rank MR diagram for HHG that
satisfy the  conditions
of theorem 3.4.

By the same argument that was used to prove the analogous claim in theorem 3.3,
if $Out(G)$ is infinite, then at least one of the $m$-collections in the higher rank MR diagram   
of hybrid resolutions contains a hybrid resolution in which at least one of the   resolutions in
the hybrid resolution (the one that is associated with the top level, or one of the resolutions that are associated with the domains (set) stabilizers)
has
 at least two levels. 

\line{\hss$\qed$}

\vglue 1.5pc
\centerline{\bf{\S4. Homomorphisms and finitely generated groups}}
\medskip

In the first section we constructed weakly acylindrical actions of RAAGs on  simplicial trees, that can serve to analyze their
automorphisms. In the second section we analyzed automorphisms of groups that act on products of hyperbolic spaces, and constructed a 
higher rank Makanin-Razborov diagram that encodes the automorphisms of these groups. In the third section we generalized the construction
to colorable HHGs with weakly acylindrical actions on their domains.

In this section, we apply the constructions that appear in sections 2 and 3, to associate a higher rank Makanin-Razborov diagram with the set of
homomorphisms from a f.g.\ group into a colorable HHG. We further associate a higher rank MR diagram with every f.g.\  group, that encodes
its actions on the class of uniformly colorable HHSs with uniformly weakly acylindrical actions of their isometry groups on their domains.

\vglue 1.5pc
\proclaim{Theorem 4.1 (cf. Theorem 3.4)} Let $G$ be a  colorable HHG that acts properly and cocompactly on an HHS $X$, and suppose that the action satisfies
the assumptions of theorem 3.4. In particular, that the action of $G$ has $m$ orbits of domains of $X$, and each orbit
is pairwise transverse, and that the action of $G$ satisfies the weakly acylindrical assumptions that are listed in the statement of theorem 3.4.

Let $\Gamma$ be a f.g.\ group. With the set of homomorphisms $Hom(\Gamma,G)$ it is possible to associate a finite set of $m$-collections of hybrid resolutions, where each
hybrid resolution has two parts, precisely as in the higher rank diagram for $Aut(G)$ in theorem 3.4. 

Every homomorphism from $\Gamma$ to $G$ factors through at least one of the $m$-collections of hybrid resolutions (see definition 2.6 for a homomorphism that factors through
an $m$-collection of (hybrid) resolutions). If $Hom(\Gamma,G)$ has infinitely
many non-conjugate homomorphisms (non-conjugate in $G$), then 
at least one of the $m$-collections  of hybrid  resolutions contains a hybrid resolution in which at least one of the   resolutions that the hybrid resolution
 is composed from  has
 at least two levels.
\endproclaim

\nfp 
The proof is identical to the proof of theorem 3.4. Indeed, the proof of theorem 3.4 did not use the automorphism assumption, nor that the
domain was $H$ itself.

\line{\hss$\qed$}

The next theorem is what we view as an analogue of Thurston's bounded image theorem on the geometric structure of discrete faithful
representations of a f.g.\ group into a rank 1 Lie group and its connection to the JSJ decomposition (e.g., [Mo], [Se3] or [Ka]), for uniformly
 weakly acylindrical actions of f.g.\ groups on 
uniformly colorable HHS.

\proclaim{Theorem 4.2} Let $\Gamma$ be a f.g.\ group, and let  $m$ be a positive integer. 
There exists a higher rank Makanin-Razborov diagram, similar to the one that
was constructed in theorem 3.4 (i.e., with finitely many $m$-collections of hybrid resolutions), that is constructed from all the homomorphisms of  $\Gamma$ into 
all the groups $G$ that act an HHS $X$, where the family of HHS $X$ and the actions of $G$ on each of its members satisfy the following conditions:
\roster
\item"{(1)}" there are exactly $m$ orbits of domains of an HHS $X$ under the action of $G$. $X$ is colorable w.r.t. the action of $G$.

\item"{(2)}" the family of HHS $X$ have uniform HHS structural constants.

\item"{(3)}" $G$ acts on the $m$ projection complexes that are associated with its action on $X$ uniformly weakly acylindrically (i.e., there are global 
 acylindricity
constants $r,b$ for all the family of actions on the associated projection complexes).

\item"{(4)}" the domain (set) stabilizers in $G$ act uniformly weakly acylindrically on the domain in $X$ that they stabilize.
\endroster  
\endproclaim

\nfp Identical to the proof of theorem 3.4.

\line{\hss$\qed$}


\vskip 2em

\end